\documentclass[a4]{article} 
\usepackage{amssymb} 
\usepackage{latexsym} 
\usepackage{psfig} 
\usepackage{graphicx}
\usepackage{psfrag}
\newcommand{\newtext}{} 
\newcommand{\nottopapertext}{} 

\newcommand{\R}{{\mathbb R}}

\newcommand{\bt}{{\bf t}} 
\newcommand{\bn}{{\bf n}} 
\newcommand{\ud}{{\underline{\delta}}} 
\newcommand{\Hnull}{{\stackrel{\circ}H}{}}  
\newcommand{\Hbnull}{{\stackrel{\circ}{\bf H}{}}}  
\newcommand{\Cnull}{{{\stackrel{\circ}C}{}^\infty}}

\def\hat{\widehat}
\def\tilde{\widetilde}
\def \bfo {\begin {eqnarray*} }
\def \efo {\end {eqnarray*} }
\def \ba {\begin {eqnarray*} }
\def \ea {\end {eqnarray*} }
\def \beq {\begin {eqnarray}}
\def \eeq {\end {eqnarray}}

\def \diam {\hbox{diam }}

\def \e {\varepsilon}
\def \p {\partial}

\def\M{{\mathcal M}}
\def\F{{\mathcal F}}

\def\hf{{\hat f}}
\def\hh{{\hat h}}
\def\ud{{\underline \delta}}

\newtheorem{definition}{Definition}[section] 
\newtheorem{theorem}[definition]{Theorem} 
\newtheorem{lemma}[definition]{Lemma} 
 
\newtheorem{proposition}[definition]{Proposition} 
\newtheorem{corollary}[definition]{Corollary} 
 
\begin{document} 
\title{Maxwell's Equations with Scalar Impedance: 
Direct and Inverse Problems} 
\author{Yaroslav Kurylev 
\thanks{Loughborough University} 
\and Matti Lassas\thanks{University of 
Helsinki} \and Erkki Somersalo\thanks 
{Helsinki University of Technology}} 
\date{28 Oct. 2002}
\maketitle

{\bf Abstract:} The article deals with electrodynamics in the
presence of anisotropic materials having scalar wave impedance.
Maxwell's equations written for differential forms over a 
3-manifold are analysed. The system is extended to a Dirac type first order
elliptic system on the Grassmannian bundle over the manifold. 
The second part of the article deals with the dynamical inverse
boundary value problem of determining the electromagnetic material
parameters from boundary measurements. By using the boundary control
method, it is proved that 
the dynamical boundary data determines
the electromagnetic travel time metric as well as the scalar
wave impedance on the manifold. This invariant result leads also to a complete
characterization of the non-uniqueness of the corresponding
inverse problem in bounded domains of $\R^3$. 

\section*{Introduction}\label{introduction}

Classically, the laws of electromagnetism expressed by
Maxwell's equations are written for
vector fields representing the electric and magnetic fields.
However, it is possible to rephrase these equations
in terms of differential forms. It turns out that this
alternative formulation has several advantages both
from the theoretical and practical point of view.
First, the formulation of electromagnetics with
differential forms reflect the way in which the
fields are actually observed. For instance, flux
quantities are expressed as 2--forms while field
quantities that correspond to forces are naturally
written as 1--forms. 
This point of view has been adopted in modern physics
at least when fields in free space are dealt with, see
\cite{Frankel}.
Furthermore, the formulation distinguishes the
topological properties of 
 the electromagnetic media
 from
those that depend on geometry. It is understood that geometry is
related to the properties of the material where the waves propagate.
The distinction between non-geometric and geometric
properties has consequences also to the numerical
treatment of the equations by so called Whitney forms.
An extensive treatment of this topic can be found in
\cite{bossavit1}, \cite{bossavit2}. 
For the original reference concerning Whitney elements
see \cite{whitney}.

The present work is divided in two parts.
In the first part, we pursue further the invariant formulation of
Maxwell's equations to model the wave propagation
in certain anisotropic materials. More
precisely, we consider anisotropic materials with scalar wave
impedance. Physically, scalar wave impedance is tantamount to
a single propagation speed of waves with different polarization.
The invariant approach leads us to formulate Maxwell's equations on
3-manifolds as a first order Dirac type system. From the operator
theoretic point of view, this formulation is based on an
elliptization procedure by extending Maxwell's equations to a 
Grassmannian bundle over the manifold. This is a generalization of
the elliptization of Birman 
and Solomyak
and Picard (see\cite{Birman},\cite{picard}).

In the second part of the work, we consider
the inverse boundary value problem 
for
Maxwell's equations.
In terms of physics,
the
 goal is to determine material parameter tensors, 
electric permittivity $\epsilon$ and magnetic permeability $\mu$,
in a bounded
domain from field observations at the boundary of that
domain. As it is already well established, for anisotropic
inverse problems it is natural to consider the problem in
two parts. First, we consider the invariant problem on
a Riemannian manifold, where we recover the travel time
metric and the wave impedance on the manifold. As a second
step, we consider the consequences of the invariant result
when the manifold is imbedded to $\R^3$.

Although inverse problems in electrodynamics have a great
significance in physics and applications, results concerning
the multidimensional inverse problems are relatively recent.
One-dimensional results have existed starting from the 30'ies,
see e.g. \cite{langer}, \cite{schlichter}. 
The first breakthrough in multidimensional
inverse problems for electrodynamics was based on the
use of complex geometrical optics \cite{SIC}, \cite{CP}, \cite{OPS},
\cite{OS}.  In these papers, the inverse problem
of recovering the scalar material parameters from complete
fixed frequency boundary data was solved even in the non-selfadjoint
case, i.e., 
in the
presence of electric conductivity.
These works were based on ideas previously developed in references
\cite{SyU},\cite{Na1},\cite{Na2} to solve the 
scalar Calder\'{o}n problem, that obtained its present formulation
in \cite{Cl}.

In the dynamical case, 
a method to solve an isotropic inverse boundary problem based on ideas
of integral geometry  is developed in \cite{rom}. The method, however,
is confined to the case of a geodesically simple manifolds and, at the moment,
is limited to finding some combinations of material parameters, including
electric conductivity.
An alternative method to tackle the
inverse boundary value problem is the boundary control (BC) method,
originated in \cite{Be1}. Later, this method was developed for
the Laplacian on Riemannian manifolds \cite{BeKu3} and for 
anisotropic 
self-adjoint \cite{Ku1} -- \cite{Ku3} and certain 
non-selfadjoint inverse problems \cite{KL2}. The first
application of the BC method to electrodynamics was done 
in \cite{Be4}, \cite{BeIsPSh}.
The authors of these articles     
show that, when the material parameters $\epsilon$ and
$\mu$ are real scalars or alternatively when $\epsilon=\mu$, the boundary
data determines the wave speed in the vicinity of the boundary.
These works employed the Hodge-Weyl decomposition in the domain
of influence near the boundary. The real obstruction for this technique
is that, as time grows, the domain of influence can become 
non-smooth and the topology may be highly involved. For these
reasons, our paper is based on different ideas.

In this article, there are essentially two new leading ideas. First,
we characterize the subspaces controlled from the boundary by
duality, thus avoiding the difficulties arising from the complicated
topology of the domain of influence. The second idea is to develop
a method 
of  waves focusing at
a single point of the manifold. This
enables us to recover pointwise values of the waves on the manifold.
The geometric techniques of the paper are presented 
in \cite{Ku5} and
 the book
\cite{KKL}.

The main results of this paper can be summarized as follows.
\begin{enumerate}
\item The knowledge of the complete dynamical boundary data 
over a 
sufficiently large
finite period of time determines uniquely
the compact manifold endowed with the electromagnetic travel time metric
as well as the scalar wave impedance (Theorem 4.1).
\item For the corresponding anisotropic
inverse boundary value problem with scalar wave impedance for
bounded domains in $\R^3$, the non-uniqueness 
is completely characterized by describing
the class of possible transformations between material tensors
that are indistinguishable from the boundary (Theorem 11.1). 
\end{enumerate} 
To the best knowledge of the authors, no global uniqueness
results for inverse problems for systems with anisotropic 
coefficients have been previously known.

\bigskip
{\bf Acknowledgements:} We would like to give our warmest thanks to
professor Alexander Katchalov for numerous useful discussions.
His lectures on non-stationary Gaussian beams 
\cite{Ka1}, \cite{Ka2} at Helsinki
University of Technology were paramount for our understanding
of the subject.
This work was accomplished during several
visits of the authors at each others' home institutions.
We wish to thank Helsinki University of Technology, Loughborough
University and University of Helsinki for their kind hospitality
and financial support.
Furthermore, the financial support of the Academy of Finland and
Royal Society is acknowledged.

\section{Maxwell's equations for forms}\label{maxwell for forms} 

 In this chapter we derive 
an
invariant form
 for Maxwell equations,
consider initial boundary value problem for them and show 
how energy of fields can be found using boundary measurements.

We start with Maxwell equations in domain $\Omega\subset \R^3$
equipped with the standard Euclidean structure. Since our objective
is to write Maxwell equations in 
 an 
invariant form, we generalize
the setting in very beginning and instead of domain $\Omega$
consider manifolds. 
 
Let $(M,g_0)$ be a 
connected, oriented Riemannian 
3-manifold possibly with a boundary $\partial M\neq 
\emptyset$. We assume that all objects in this paper are $C^\infty$--smooth.  
Consider Maxwell's equations on $M$, 
\begin{eqnarray} 
 {\rm curl}\,E &=& - B_t,\mbox{ (Maxwell--Faraday)},\label{MF vector}\\ 
\noalign{\vskip4pt} 
 {\rm curl }\,H &=& \phantom{-}D_t,\mbox{ (Maxwell--Amp\`{e}re)}, 
\label{MA vector} 
\end{eqnarray} 
where $E$ and $H$ are the electric and magnetic fields, 
and $B$ and $D$ are the magnetic flux density 
and electric displacement, assumed for the time being 
to be smooth mappings $M\times\R\to TM$. Here  
$TM$ denotes the tangent bundle over $M$.  
The curl operator as well as divergence appearing later
will be defined invariantly in formula (\ref{A 23}) below.
The sub-index $t$ in the equations  
(\ref{MF vector})--(\ref{MA vector}) 
denotes differentiation with respect to time. 
We denote 
the collection of these vector fields as $\Gamma(M\times\R)$.  
At this point, we do not specify the 
initial and boundary values. To avoid non-physical 
static solutions, the above equations are  
augmented with the conditions 
\begin{equation}\label{div eqs} 
 {\rm div}B=0,\quad {\rm div}D=0. 
\end{equation} 
 
Furthermore, the fields $E$ and $D$, and similarly 
the fields $H$ and $B$ are interrelated through the 
constitutive relations. In anisotropic and non-dispersive 
medium, the constitutive relations assume the simple 
form 
\begin{equation}\label{constitutive} 
 D=\epsilon E,\quad B=\mu H, 
\end{equation} 
where $\epsilon,\mu$ are smooth and strictly positive 
definite tensor fields of type $(1,1)$ on $M$. 
Our aim is to write the above equations using differential 
forms. 
 
Given the metric $g_0$, we can associate in a canonical way 
a differential 1--form to correspond each vector field. Let us denote 
by $\wedge^k T^*M$ 
the $k$:th exterior power of the cotangent bundle. 
We define the mapping 
\[ 
 TM\to T^*M,\quad X\mapsto X^\flat 
\] 
through the formula $g_0(X,Y) = X^\flat(Y)$.  
This mapping is one-to-one and it has the following 
well-known properties (See e.g. \cite{Sc}): For 
a scalar field $u\in C^\infty(M)$, 
$ 
({\rm grad}\,u)^\flat = du, 
$ 
where $d$ is the exterior differential 
and for a vector field $X\in\Gamma(M)$, we have 
\beq \label{A 23}
({\rm curl}\,X)^\flat = *_0 dX^\flat,\quad 
{\rm div}\,X = -\delta_0 X^\flat, 
\eeq
where $*_0$ denotes the Hodge--$*$ operator 
with respect to the metric $g_0$, 
\[ 
 *_0:\wedge^k T^*M\to\wedge^{3-k}T^*M, 
\] 
and $\delta_0$ denotes the codifferential 
\footnote{Cf. with $\delta_0=(-1)^{nk+n+1}*_0d*_0$ for 
Riemannian $n$--manifolds}, 
\[ 
 \delta_0 = (-1)^k*_0 d *_0: 
\Omega^k M\rightarrow \Omega^{k-1} M. 
\] 
Here, $\Omega^k M$ denotes the smooth sections 
$M\to\wedge^k T^* M$,
i.e. differential $k-$forms.
Applying now the operator $\flat$ on Maxwell's equations 
(\ref{MF vector})--(\ref{MA vector}) yields 
\[ 
dE^{\flat} = - *_0B^{\flat}_t,\quad 
 dH^{\flat} = *_0D^{\flat}_t, 
\] 
where we used the identity $*_0*_0 = {\rm id}$ valid in 
3--geometry\footnote{For Riemannian $n$--manifold, 
we have in general $*_0*_0=(-1)^{k(n-k)}$}. 
The divergence equations (\ref{div eqs}) read 
\[ 
 \delta_0 D^{\flat}=0,\quad \delta_0 B^{\flat}=0. 
\] 
 
Consider now the constitutive relations 
(\ref{constitutive}). Starting with the equation 
$D=\epsilon E$, we pose the following question: 
Is it possible to find a {\em metric} 
$g_\epsilon$ such that the Hodge-$*$ operator 
with respect to this metric, denoted by $*_\epsilon$, 
would satisfy the identity 
\[ 
 *_0 D^\flat =*_0(\epsilon E)^\flat = 
 *_\epsilon E^\flat? 
\] 
Assume that such a metric $g_\epsilon$ exists. 
By writing out the above formula 
in given 
{\nottopapertext local coordinates $(x^1,x^2,x^3)$}
and recalling 
the definition of the Hodge-$*$ operator, the left side 
yields 
\begin{eqnarray*} 
 *_0(\epsilon E)^\flat &=& 
 \sqrt{{\rm g}_0} g_0^{ij}e_{jpq}g_{0,ij} 
\epsilon^j_kE^k dx^p\wedge dx^q \\ 
\noalign{\vskip6pt} 
&=&\sqrt{{\rm g}_0}e_{jpq}\epsilon^j_k E^k 
 dx^p\wedge dx^q, 
\end{eqnarray*} 
where $e$ is the totally antisymmetric permutation 
index and ${\rm g}_0={\rm det}(g_{0,ij})$. 
Likewise, the right side reads 
\[ 
 *_\epsilon E^\flat = \sqrt{{\rm g}_\epsilon} 
 g_\epsilon^{ij}e_{jpq}g_{0,ik} E^k dx^p\wedge 
 dx^q, 
\] 
so evidently the desired equality ensues if 
we set 
\[ 
 \sqrt{{\rm g}_\epsilon} 
 g_\epsilon^{ij}g_{0,ik} =\sqrt{{\rm g}_0} 
 \epsilon_j^k. 
\] 
By taking determinants of both sides we find 
that 
\[ 
 \sqrt{{\rm g}_\epsilon}=\sqrt{{\rm g}_0} 
 {\rm det}(\epsilon). 
\] 
Thus
we see that the appropriate form for the metric tensor in the 
contravariant form is 
\beq\label{connection e and g}
 g_\epsilon^{ij} = \frac 1{{\rm det}(\epsilon)} 
g_0^{ik}\epsilon^j_k. 
\eeq
In the same fashion, we find a metric $g_\mu$ 
such that  
\[ 
 *_0 B^\flat =*_0(\mu H)^\flat = 
 *_\mu H^\flat. 
\] 
In general, the metrics $g_\mu$ and $g_\epsilon$ 
can be very different from each other. In this 
article, we consider a particular case. Indeed, 
assume that the material has a {\em scalar 
wave impedance}, i.e., the tensors $\epsilon$ 
and $\mu$ satisfy 
\[ 
 \mu =\alpha^2 \epsilon, 
\] 
where the wave impedance, $\alpha =\alpha(x)$,
is a smooth function on $M$. Now we define 
two families of 1-- and 2--forms on $M$ as 
follows. We set 
\[ 
 \omega^1 = E^\flat,\quad \omega^2 = *_0B^\flat. 
\] 
Similarly, we define 
\beq\label{nu field 1}
 \nu^1 = \alpha H^\flat,\quad \nu^2 =*_0\alpha D^\flat. 
\eeq
Observe that the wave impedance scaling renders 
$\omega^1$ and $\eta^1$ to have the same physical dimensions, 
and the same holds for the 2--form. 
Now it is a straightforward matter to check that the 
constitutive relations assume the form 
\[ 
 \nu^2 =\alpha*_\epsilon\omega^1,\quad 
 \omega^2=\frac 1\alpha *_\mu \nu^1. 
\] 
We can make these equations even more symmetric by 
proper scaling of the metrics. Indeed, since 
$
\alpha^{-1}\mu =\alpha\epsilon, 
$ 
we have a new metric $g$ that is defined as 
\[ 
 g^{ij}=g_{\alpha\epsilon}^{ij}=g_{\alpha^{-1}\mu}^{ij}. 
\] 
We have, by direct substitution that 
\begin{equation}\label{travel time metric} 
  g^{ij}=\frac 1{\alpha^2}g_{\epsilon}^{ij} 
 =\alpha^2 g_{\mu}^{ij}. 
\end{equation} 
This new metric will be called the {\em travel time 
metric} in the sequel. 
 
Assume that 
\[ 
 *:\wedge^j T^* M\to \wedge^{3-j}T^* M 
\] 
denotes the Hodge--$*$ operator with respect to 
some
metric $\hat{g}$. 
{\nottopapertext
 If we perform a scaling of the metric as 
\[ 
 \hat{g}^{ij}\rightarrow \tilde g^{ij}=r^2\hat{g}^{ij}, 
\] 
the corresponding Hodge operator is scaled as 
\[ 
 *\rightarrow \tilde * = r^{2j-3}*. 
\] 
Therefore, if we denote by $*$ the Hodge--$*$ operator 
with respect to the travel time metric, we have 
}
\[ 
 * = \alpha *_\epsilon = \frac 1\alpha *_\mu : 
 \wedge^1 T^* M\to \wedge^2 T^* M. 
\] 
But this means simply that, in terms of the travel 
time metric, we have 
\begin{equation}\label{hodge} 
 \nu^2 =*\omega^1,\quad \omega^2=*\nu^1. 
\end{equation} 
 
Consider now Maxwell's equations for these forms. 
After eliminating the $\nu$--forms using the  
constitutive equations (\ref{hodge}), 
Maxwell--Faraday and Maxwell--Amp\`{e}re equations 
assume the form 
\begin{equation}\label{apu1} 
 d\omega^1=-\omega_t^2,\quad 
 \delta_\alpha\omega^2 = \omega_t^1, \quad  
 \delta_\alpha =
 (-1)^k*\alpha d  \frac 1\alpha*: 
\Omega^k M\rightarrow \Omega^{k-1} M
\end{equation} 
and the divergence equations (\ref{div eqs}) read 
\begin{equation}\label{apu2} 
 d\omega^2 =0,\quad \delta_\alpha\omega^1=0. 
\end{equation} 
{\newtext In the sequel, we call equations (\ref{apu1})
and (\ref{apu2}) {\it Maxwell's equations}.}

It turns out to be useful to define auxiliary 
forms that vanish in the electromagnetic theory. 
Let us introduce the auxiliary forms 
$\omega^0$ and $\omega^3$ via the formulas 
\[ 
 \omega_t^0 = \delta_\alpha\omega^1,\quad 
 -\omega^3_t=d\omega^2. 
\] 
Furthermore, we define the corresponding 
$\nu$--forms as 
\beq\label{nu field 2}
 \nu^0=*\omega^3,\quad \nu^3=*\omega^0. 
\eeq
Since these auxiliary 
forms
 are all vanishing, 
we may modify the equations (\ref{apu1}) to have 
\begin{equation}\label{apu1 bis} 
 d\omega^1 -\delta_\alpha\omega^3=-\omega_t^2, 
\quad  
 d\omega_0 -\delta_\alpha\omega^2  = -\omega_t^1. 
\end{equation} 
Putting the obtained equations together in a matrix 
form, we arrive at the equation 
\begin{equation}\label{complete} 
 \omega_t + {\cal M}\omega=0, 
\end{equation} 
where  
\[ 
 \omega = (\omega^0,\omega^1,\omega^2,\omega^3) 
\] 
and the operator ${\cal M}$ (without defining its domain at this 
point, i.e., defined as a differential expression) 
is given as 
\begin{equation}\label{M} 
{\cal M} = \left(\begin{array}{cccc} 
0 &-\delta_{\alpha} &0       &0       \\ 
d & 0      &-\delta_{\alpha} &0       \\ 
0 & d      &0       &-\delta_{\alpha} \\ 
0 & 0      &d       &0 
\end{array}\right). 
\end{equation} 
The equation (\ref{complete}) is called {\em 
the complete Maxwell system}. 
In the next section, we treat more systematically 
this operator.  
 
{\bf Remark 1.} The operator ${\mathcal M}$ has the 
property 
\[ 
 {\mathcal M}^2 = -{\rm diag}(\Delta_\alpha^0, 
\Delta_\alpha^1,\Delta_\alpha^2,\Delta_\alpha^3) 
 =-{\bf \Delta}_\alpha, 
\] 
where the operator $\Delta_\alpha^k$ acting on 
$k$--forms is 
\[ 
 \Delta_\alpha^k = d\delta_\alpha + \delta_\alpha d 
 =\Delta^{k}_g + Q(x,D), 
\] 
with $\Delta_g^k$ denoting the Laplace-Beltrami operator 
on 
$k$--forms
with respect to the travel time metric and $Q(x,D)$ 
being a first order perturbation. 
Hence, if $\omega$ satisfies the equation 
(\ref{complete}), we have 
\[ 
 (\partial_t^2 + {\bf \Delta}_\alpha)\omega 
 =(\partial_t-{\mathcal M})(\partial_t 
+{\mathcal M})\omega =0. 
\] 
In particular, we observe that the assumption that 
the impedance is scalar implies a unique propagation 
speed for the system. 
\smallskip

{\bf Remark 2.} Denote by $\Omega M=\oplus_{k=0}^3 \Omega^k M$
the Grassmannian algebra of differential forms, 
where $\Omega^k M$ are the differential $k$-forms.
 Then the operator $\M$ in formula 
(\ref{M}) can be also considered
as a Dirac operator $d-\delta_\alpha:\Omega M\to \Omega M$.
\smallskip

Before leaving this section, let us 
briefly consider the energy integrals in terms of 
the differential forms. In terms of the vector fields, the 
energy of the electric field at a given moment $t$ 
is obtained as the integral 
\[ 
 {\mathcal E}(E) = \int_M \epsilon E\cdot E dV 
 =\int_M g_0(E,D)dV = \int_M E^\flat\wedge *_0D^\flat
\] 
where $dV$ is volume form of $(M,g_0)$.
By plugging in the defined forms we arrive at 
\[ 
 {\mathcal E}(E) =\int_M\frac 1\alpha \omega^1\wedge 
 *\omega^1. 
\] 
In the same fashion, we find that the energy of the 
magnetic field reads 
\[ 
 {\mathcal E}(B) =\int_M\frac 1\alpha \omega^2\wedge 
 *\omega^2. 
\] 
These formulas serve as a motivation for our definition 
of the inner product in the following section. 
 
\subsection{Maxwell operator} 
 
In this section we establish a number of notational 
conventions and definitions concerning the differential 
forms used in this work. 
 
We define the $L^2$--inner products for $k$--forms 
in $\Omega^k M$ as 
\[ 
 (\omega^k,\eta^k)_{L^2} = \int_M \frac 1\alpha 
 \omega^k\wedge * \eta^k,\quad 
 \omega^k,\;\eta^k\in\Omega^k M. 
\] 
Further, we denote by 
 $L^2(\Omega^k M)$ the completion of 
$\Omega^k M$ with respect to the norm 
defined by the above inner products. 
We also define 
\[ 
 {\bf L}^2(M) =  
 L^2(\Omega^0 M)\times  
L^2(\Omega^1 M)\times  
L^2(\Omega^2 M)\times L^2(\Omega^3 M). 
\] 
Similarly, we define Sobolev spaces ${\bf H}^s(M),  \, s \in \Bbb{R}$,
\[ 
{\bf H}^s(M) = H^s(\Omega^0M)\times  
H^s(\Omega^1 M)\times  
H^s(\Omega^2 M)\times H^s(\Omega^3 M),
\]
\[ 
{\bf H}_0^s(M) = H^s_0(\Omega^0M)\times  
H^s_0(\Omega^1 M)\times  
H^s_0(\Omega^2 M)\times H^s_0(\Omega^3 M),
\]
where $H^s(\Omega^kM)$ are Sobolev spaces of $k-$ forms. 
At last, $H^s_0(\Omega^k M)$ is the closure in $H^s(\Omega^k M)$
of $\Omega^k M^{{\rm int}}$, i.e. the subspace of $\Omega^k M$
of $k-$ forms which vanish near $\p M$. 
 
The domain of the exterior derivative $d$ in the 
$L^2$--space of $k$--forms is 
\[ 
 H(d,\Omega^k M) = \left\{ \omega^k\in L^2 
 (\Omega^k M)\mid 
 d\omega^k \in L^2(\Omega^{k+1} M)\right\}. 
\] 
Similarly, we set 
\[ 
 H(\delta_\alpha,\Omega^k M) = \left\{ \omega^k 
 \in L^2(\Omega^k M)\mid 
 \delta_\alpha\omega^k \in L^2(\Omega^{k-1} M)\right\}, 
\] 
where $\delta_\alpha$ is the weak extension of the 
operator 
$
\delta_\alpha:\Omega^k M\to\Omega^{k-1}M. 
$ 
In the sequel, we shall drop the sub-index $\alpha$ 
from the codifferential. 
 
The codifferentiation $\delta$ is adjoint to the 
exterior derivative in the sense that for $C^\infty 
_0$--forms on $M$, 
\[ 
(d\omega^k,\eta^{k+1})_{L^2} = (\omega^k,\delta\eta^{k+1})_{L^2}. 
\] 
To extend the adjoint formula for less regular forms, 
let us first fix some notations. For $\omega^k 
\in\Omega^k M$, we define the {\em tangential} 
and {\em normal} boundary data at $\partial M$ as 
\begin{eqnarray*} 
\bt\omega^k = i^*\omega^k,\quad 
\bn\omega^k = i^*(\frac 1\alpha*\omega^k), 
\end{eqnarray*} 
respectively, where $i^*:\Omega^k M 
\to\Omega^k\partial M$ is the 
pull-back of the natural imbedding 
$i:\partial M\to M$.  
 Sometimes, we denote $\bn=\bn_\alpha$ to indicate a particular 
 choice
$\alpha$.
With these notations, let us write 
\[ 
  \int_{\partial M}i^*\omega^k 
\wedge i^*(\frac 1\alpha 
* \eta^{k+1})=\langle\bt\omega^k, 
\bn\eta^{k+1}\rangle. 
\] 
We add here a small {\em caveat} that the 
above formula does not define an inner 
product as $\omega^k$ and $\eta^{k+1}$ are differential 
forms of different order. 
For $\omega^k\in\Omega^kM$ and $\eta^{k+1}\in 
\Omega^{k+1}M$, 
the Stokes formula for forms can be written as 
\begin{equation}\label{stokes1} 
 (d\omega^k,\eta^{k+1})_{L^2}- 
(\omega^k,\delta 
 \eta^{k+1})_{L^2} =  
 \langle\bt \omega^k, 
\bn \eta^{k+1}\rangle. 
\end{equation} 
{\nottopapertext 
This formula allows the extension of the boundary 
trace operators $\bt$ and $\bn$ to $H(d,\Omega^k M)$ 
and $H(\delta,\Omega^k M)$, respectively. Indeed, if 
$\omega^k\in H^1(\Omega^k M)$,  
then 
 $\bt\omega^k\in  
H^{1/2}(\Omega^k\partial M)$ 
and, by formula (\ref{stokes1}), we may extend 
\[ 
 \bt:H(d,\Omega^k M)\to H^{-1/2}(\Omega^k\partial M). 
\] 
In the same way, equation (\ref{stokes1}) gives us 
the natural extension 
\[ 
 \bn: H(\delta,\Omega^{k+1} M) 
= H^{-1/2}(\Omega^{2-k}\partial M), 
\] 
In fact, a stronger result holds. 
} 
\begin{proposition}\label{paquet} 
The operators $\bt$ and $\bn$  
can be extended to continuous surjective maps 
\begin{eqnarray*} 
\bt: H(d,\Omega^k M) &\to&  
H^{-1/2}(d,\Omega^k\partial M), \\ 
\bn: H(\delta,\Omega^{k+1} M) &\to&  
H^{-1/2}(d,\Omega^{2-k}\partial M), \\ 
\end{eqnarray*} 
where the space $H^{-1/2}(d,\Omega^k\partial M)$ is the 
space of $k$-forms $\omega^k$  
on $\partial M$ satisfying 
\[ 
 \omega^k\in H^{-1/2}(\Omega^k \partial M),\quad 
 d\omega^k\in H^{-1/2}(\Omega^{k+1} \partial M). 
\] 
\end{proposition} 
 
This result is due to Paquet \cite{paquet}. 
 
The formula (\ref{stokes1}) can be used 
also to define function spaces with vanishing 
boundary data. Indeed, let us define 
\begin{eqnarray*} 
 \Hnull (d,\Omega^k M) =\{ \omega^k\in 
H(d,\Omega^k M) &\mid& (d\omega^k,\eta^{k+1})_{L^2} 
=(\omega^k,\delta\eta^{k+1})_{L^2} \\ 
 & &\mbox{for all 
$\eta^{k+1}\in H(\delta,\Omega^{k+1}M)$}\},\\
 \Hnull(\delta,\Omega^{k+1} M)_ =\{ \eta^{k+1}\in 
H(\delta,\Omega^{k+1} M)&\mid& (d\omega^k,\eta^{k+1})_{L^2} 
=(\omega^k,\delta\eta^{k+1})_{L^2} \\ 
 & &\mbox{for all 
$\omega^k\in H(d,\Omega^kM)$}\}. 
\end{eqnarray*} 
It is not hard to see that indeed 
\[ 
 \Hnull(d,\Omega^k M)=\bt^{-1}\{0\}, 
\quad \Hnull(\delta,\Omega^{k+1}M)= 
\bn^{-1}\{0\}. 
\] 
 
{\nottopapertext We 
are now in the position
prove the following lemma.} 
\begin{lemma}\label{adjoints} 
The adjoint of the operator 
\[ 
 d:L^2(\Omega^k M)\supset H(d,\Omega^k M)\to L^2(\Omega^{k+1} M) 
\] 
is the operator 
\[ 
 \delta:L^2(\Omega^{k+1}M) 
\supset\Hnull(\delta,\Omega^{k+1}M) 
\to L^2(\Omega^kM) 
\] 
and {\em vice versa}. Similarly, the adjoint of 
\[ 
 \delta:L^2(\Omega^{k+1}M) 
\supset H(\delta,\Omega^{k+1}M) 
\to L^2(\Omega^kM) 
\] 
is the operator 
\[ 
 d:L^2(\Omega^k M)\supset \Hnull(d,\Omega^k M)\to L^2(\Omega^{k+1} M) 
\] 
\end{lemma} 

{\nottopapertext 
{\em Proof:} We prove only the first claim, the other 
having a similar proof. 
 
Let $\eta^{k+1}\in{\cal D}(d^*)$,  
where $d^*$ denotes the 
adjoint of $d$. By definition, there exists 
$\vartheta^k\in L^2(\Omega^k M)$ such that 
\[ 
 (d\omega^k,\eta^{k+1})_{L^2} = 
(\omega^k,\vartheta^k)_{L^2} 
\] 
for all $\omega^k\in H(d,\Omega^k M)$. In particular, 
if $\omega^k\in\Omega^k M^{{\rm int}}$, we see that, in the weak sense, 
\[ 
 (d\omega^k,\eta^{k+1})_{L^2} =(\omega^k, 
\delta\eta^{k+1})_{L^2} = 
(\omega^k,\vartheta^k)_{L^2}, 
\] 
i.e., $\delta\eta^{k+1}=\vartheta^k 
\in L^2(\Omega^k M)$. Thus, $\eta^{k+1} 
\in H(\delta,\Omega^{k+1} M)$, and the 
claim follows now since we have 
\[ 
 (d\omega^k,\eta^{k+1})_{L^2} = 
(\omega^k,\delta\eta^{k+1})_{L^2} 
\] 
for all $\omega^k\in H(d,\Omega^k M)$, 
i.e., $\delta = d^*$. 
\hfill$\Box$ 
}
 
In the sequel, we will write for brevity 
$H(d)=H(d,\Omega^k M)$, etc. when there is no 
risk of confusion concerning the order of the 
forms. 
 
For later reference, let us point out that the 
Stokes formula for the complete Maxwell system 
can be written compactly as 
\begin{equation}\label{stokes for system} 
 (\eta,{\cal M}\omega)_{{\bf L}^2} +  ({\cal M} \eta, 
\omega)_{{\bf L}^2} = 
\langle\bt\omega,\bn\eta\rangle +  
\langle\bt\eta,\bn\omega\rangle, 
\end{equation} 
where $\omega \in{\bf H}$ with 
\beq
\label{H}
{\bf H} = H(d)\times[H(d)\cap 
H(\delta)]\times[H(d)\cap 
H(\delta)]\times H(\delta) 
\eeq 
and $\eta \in {\bf H}^1(M)$ and
we  use the notations 
\begin{eqnarray*} 
\bt\omega = (\bt\omega^0,\bt\omega^1,\bt\omega^2)\
\quad
\bn\omega = (\bn \omega^3,\bn \omega^2,\bn 
 \omega^1), 
\end{eqnarray*} 
 and, naturally, 
\[ 
\langle\bt\omega,\bn\eta\rangle = 
\langle\bt\omega^0, 
\bn\eta^1\rangle + \langle\bt\omega^1, 
\bn\eta^2\rangle +\langle\bt\omega^2, 
\bn\eta^3\rangle. 
\] 
 
With these notations, we give the following definition 
of the Maxwell operators with electric and magnetic 
boundary conditions, respectively. 
 
 
\begin{definition}\label{d. 2} The Maxwell operator with the 
electric boundary condition, denoted by  
\[ 
 {\cal M}_{\rm e}:{\cal D}({\cal M}_{\rm e})\to 
{\bf L}^2(M), 
\] 
is defined through the differential expression 
(\ref{M}), with the domain ${\cal D}({\cal M}_{\rm e}) 
\subset {\bf L}^2(M)$ defined as 
\[ 
 {\cal D}({\cal M}_{\rm e})=\Hbnull_\bt:= \Hnull(d) 
\times[\Hnull(d)\cap H(\delta)]\times 
 [\Hnull(d)\cap H(\delta)] 
\times H(\delta). 
\] 
Similarly, the Maxwell operator with the magnetic 
boundary condition, denoted 
{\bf by }
\[ 
 {\cal M}_{\rm m}:{\cal D}({\cal M}_{\rm m})\to 
{\bf L}^2(M), 
\] 
is defined through the differential expression 
(\ref{M}), with the domain ${\cal D}({\cal M}_{\rm m}) 
\subset {\bf L}^2(M)$ defined as 
\[ 
 {\cal D}({\cal M}_{\rm m})=\Hbnull_\bn:=  H(d) 
\times[H(d)\cap \Hnull(\delta)]\times 
 [H(d)\cap \Hnull(\delta)] 
\times \Hnull(\delta). 
\] 
\end{definition} 
 
Before further discussion, let us comment the boundary 
conditions in terms of physics. For vectorial  
representations of the electric and magnetic fields, 
the electric boundary condition is associated with 
electrically  
perfectly conducting boundaries, i.e., $n\times E =0$, $n\cdot B=0$, where $n$ 
is the exterior normal vector at the boundary. In 
terms of differential forms, this means simply that 
$\bt E^\flat=\bt\omega^1 = 0$ and $\bt *_0 B^\flat= 
\bt\omega^2 = 0$. On the other hand, the magnetic 
boundary conditions represent a magnetically perfectly   
conducting boundaries, i.e., $n\times H=0$, 
$n\cdot D=0$, which again in terms of forms 
reads as $\bt H^\flat = \bt(1/\alpha)\nu^1 = 0$ or  
$\bt(1/\alpha)*\omega^2 =\bn\omega^2 =0$ and 
$\bt *_0 D^\flat = \bt(1/\alpha)\nu^2=0$, or in terms of 
$\omega^1$, $\bt(1/\alpha)*\omega^1 = \bn\omega^1=0$. 
 
There is an obvious duality between these 
conditions. 
It is therefore sufficient to consider the operator 
with the electric boundary condition only. 
This observation is related to the well-known  
Maxwell {\it duality principle.} 
 
Consider the intersections of spaces appearing 
in the domains of definition in the previous definition. 
Let us denote 
\ba
& & \Hnull^1_\bt(\Omega^k M)=\{\omega^k\in H^1 
(\Omega^k M)\mid \bt\omega^k=0\},\\ 
& & \Hnull^1_\bn(\Omega^k M)=\{\omega^k\in H^1 
(\Omega^k M)\mid \bn\omega^k=0\}. 
\ea 
It is a direct consequence of Gaffney's inequality 
(see \cite{Sc}) that 
\ba
& & \Hnull(d,\Omega^k M)\cap H(\delta,\Omega^k M) 
 =\Hnull^1_\bt(\Omega^k M),\\
& & H(d,\Omega^k M)\cap \Hnull(\delta,\Omega^k M) 
 =\Hnull^1_\bn(\Omega^k M). 
\ea

The following lemma is a direct consequence of Lemma 
\ref{adjoints} and classical results
on Hodge-Weyl decomposition\cite{Sc}.  
 
\begin{lemma}\label{lem 2.4}
 The electric Maxwell operator has the following properties: 
\begin{enumerate}

\item[i.] The operator ${\cal M}_{\rm e}$ is skew-adjoint. 
\item[ii.] The operator ${\cal M}_{\rm e}$ defines an elliptic
differential operator in $M^{int}.$
\item[iii.] Ker$\,({\cal M}_{\rm e})=\{(0,\omega^1,\omega^2,
\omega^3)\in \Hbnull_\bn:\ 
d\omega^1=0,\ \delta\omega^1=0,\ 
d\omega^2=0,\ \delta\omega^2=0,\ \delta\omega^3=0\}$. 
 
\item[iv.] Ran$\,({\cal M}_{\rm e}) =
 L^2(\Omega ^0M)\times (\delta H(\delta,\Omega^2 M) +d\Hnull(d,\Omega^0 M))\times \\
 (\delta H(\delta,\Omega^2 M) +d\Hnull (d,\Omega^1 M))\times d\Hnull (d,\Omega^2 M)$.
\end{enumerate}
\end{lemma} 
 
By the skew-adjointness, it is possible to define 
weak solutions to initial-boundary-value problems 
needed later. In the sequel we denote the 
forms $\omega(x,t)$ just
by $\omega(t)$ 
when there is no danger of misunderstanding.
 
\begin{definition}\label{weak solution} 
By the {\em weak solution} to the initial boundary value 
problem 
\[ 
\omega_t+ {\mathcal M}\omega =\rho 
\in  L^1_{loc}(\R,{\bf L}^2(M)), 
\] 
\beq
\label{ibvp}
\bt\omega\big|_{\partial M\times \R}=0,\quad 
\omega(\;\cdot\;,0)=\omega_0\in{\bf L}^2, 
\eeq
we mean the form 
\[ 
 \omega(t) ={\mathcal U}(t)\omega_0+\int_0^t{\mathcal U} 
(t-s)\rho(s)ds, 
\] 
where ${\mathcal U}(t)={\rm exp}(-t{\mathcal M}_ 
{\rm e})$ is the unitary operator generated by 
${\mathcal M}_{\rm e}$. 
\end{definition} 
 
In the analogous manner, we define weak solutions with initial
data given on $t=T$,  $T\in \R$.
 Assuming $\rho \in C(\R,{\bf L}^2(M))$ and using 
 the theory of unitary groups, 
we immediately 
obtain the regularity result 
\[ 
 \omega\in C(\R,{\bf L}^2)\cap 
C^1(\R,{\bf H}').
\]

\medskip
where ${\bf H}'$ denotes the dual of ${\bf H}$. 
 
We shall need later the boundary traces of the 
weak solution. To define them, let 
$(\omega_{0n},\rho_n)\in{\mathcal D}({\mathcal M} 
_{\rm e})\times C(\R,{\mathcal D}({\mathcal M}_{\rm e}))$ 
be an approximating sequence 
of the pair 
$(\omega_0,\rho)$ 
in ${\bf L}^2 
\times C(\R,{\bf L}^2)$. 
We define 
\[ 
 \omega_n ={\mathcal U}(t)\omega_{0n}+\int_0^t{\mathcal U} 
(t-s)\rho_n(s)ds, 
\] 
whence $\omega_n\in C(\R,{\mathcal D}({\mathcal M}_{\rm e})) 
\cap C^1(\R,{\bf L}^2)$. Let 
$\varphi 
=(\varphi^0,\varphi^1,\varphi^2)$ be  
a test form, $\varphi^j\in C^\infty_0(\R, \Omega^j\partial M)$. 
Let $\eta$ be a strong solution of the initial 
boundary value problem 
\[ 
\eta_t+ {\mathcal M}\eta =0, 
\] 
\[ 
\bt\eta=\varphi,\quad \eta(\;\cdot\;,0)=0. 
\] 
We have 
\begin{eqnarray*} 
 (\eta(T),\omega_n(T))_{{\bf L}^2}&=& \int_0^T\partial_t 
 (\eta,\omega_n)dt \\ 
\noalign{\vskip6pt} 
 &=& -\int_0^T\big(({\mathcal M}\eta,\omega_n) +(\eta, 
{\mathcal M}\omega_n)\big) dt, 
\end{eqnarray*} 
and, by applying Stokes theorem, we 
deduce  
\[ 
(\eta(T),\omega_n(T))_{{\bf L}^2} = -\int_0^T\langle 
 \varphi,\bn \omega_n\rangle. 
\] 
Hence, we observe that, when going to the limit 
$n\to\infty$, the above formula 
defines 
$\bn\omega = \lim_{n\to\infty}\bn \omega_n \in 
{\mathcal D}'(\R, {\bf {\mathcal D}'}(\partial M))$,
where 
\[
 {\bf {\mathcal D}'}(\partial M) = {\mathcal D}'(\Omega^0 \partial M) \times
{\mathcal D}'(\Omega^1 \partial M) \times
{\mathcal D}'(\Omega^2 \partial M).
\]
We conclude this section with the following result. 
 
\begin{lemma}\label{weak is maxwell} 
Assume that the initial data $\omega_0$ is of the 
form $\omega_0=(0,\omega_0^1,\omega_0^2,0)$, where 
$\delta\omega_0^1 =0$, $d\omega_0^2=0$ and 
we have $\rho=0$. Then the weak solution $\omega$ 
of Definition \ref{weak solution} satisfies also  
Maxwell's equations (\ref{apu1}),  (\ref{apu2}),
i.e., $\omega^0=0$ and  
$\omega^3=0$. 
\end{lemma} 
 
{\em Proof:} 
As observed in Remark 1, 
$\omega$ and, in particular, $\omega^0$ satisfies 
the wave equation 
\[ 
 \Delta^0_{\alpha}\omega^0 + \omega^0_{tt}=0, 
\] 
in the distributional sense, along with the Dirichlet 
boundary condition $\bt\omega^0 =0$. The initial 
data for $\omega^0$ is 
\[ 
 \omega^0(0) = \omega_0^0 =0, 
\] 
and 
\[ 
 \omega^0_t(0) = \delta\omega^1\big|_{t=0} 
 = \delta\omega^1_0=0. 
\] 
Hence, we deduce that also $\omega^0=0$.
 
 Similarly, $\omega^3$ satisfies the wave equation 
with the initial data 
\[ 
 \omega^3(0) = \omega_3^0 =0, 
\] 
and 
\[ 
 \omega^3_t(0) = -d\omega^2\big|_{t=0} 
 = -d\omega^2_0=0. 
\] 
As for the boundary condition, we observe that 
\[ 
 \bt \delta\omega^3  
=\bt\omega^2_t - \bt d\omega^1  
=\partial_t\bt\omega^2 - d\bt\omega^1 =0, 
\] 
corresponding to the vanishing Neumann data for the 
function $*\omega^3$. Thus, also $\omega^3=0$. 
\hfill$\Box$ 
 
\subsection{Initial--boundary value problem} 
 
Our next goal is to consider the forward 
problem and the Cauchy data on the lateral boundary $\p M\times \R$ for 
solutions of Maxwell's equations. 
Assume that $\omega$ is a solution of the 
complete system.  
The complete Cauchy data of this solution consists of 
\[ 
 (\bt\omega(x,t),\bn\omega(x,t)) 
 ,\quad (x,t)\in\partial M\times\R_+. 
\] 
 
Assume now that $\omega$ corresponds to the 
solution of Maxwell's equations, i.e., we have 
$\omega^0=0$ and $\omega^3=0$. Consider 
 the Maxwell-Faraday equation  in
(\ref{apu1}), 
\[ 
 \omega^2_t + d\omega^1 = 0. 
\] 
By taking the tangential trace, we find that 
$
 \bt\omega^2_t= -d\bt\omega^1, 
$ 
and further,
\[ 
 \bt\omega^2(x,t)=\omega^2(0)-\int_0^t d\bt\omega^1(x,t')dt', 
\quad x\in\partial M. 
\] 
Similarly, by taking the normal trace of
the Maxwell-Amp\`{e}re 
equation in (\ref{apu1}), 
\[ 
 \omega^1_t -\delta\omega^2 =0, 
\]  
we find that 
$ \bn \omega_t^1=d\bn \omega^2, 
$
so likewise, 
\beq 
\label{2.12.1}
\bn\omega^1(x,t)= \bn \omega^1(0)+
\int_0^t d\bn\omega^2(x,t')dt', 
\quad x\in\partial M. 
\eeq     
In the sequel we shall mainly consider the case $\omega(0)=0$, when the
lateral 
Cauchy data for the original problem 
of electrodynamics is simply 
\begin{eqnarray} 
\bt\omega &=& (0,f,-\int_0^td f(t')dt'), 
\label{t-data}\\ 
\noalign{\vskip4pt} 
\bn\omega &=& (0,g,\int_0^td g(t')dt') 
\label{n-data} 
\end{eqnarray} 
where  $f$ and  
$g$ are functions of $t$ with values in $\Omega^1 \partial M$. 

The following 
theorem implies that solutions of Maxwell's equations
are solutions of the complete Maxwell system and gives sufficient
conditions for the converse result.

\begin{theorem}\label{equivalence} 
If $\omega(t)\in C(\R,{\bf H}^1)\cap 
C^1(\R,{\bf L}^2)$  
satisfies the equation 
\begin{equation}\label{M-eq} 
 \omega_t +{\cal M}\omega = 0,\quad t>0 
\end{equation} 
with vanishing initial data $\omega(0)=0$,
and $\omega^0(t)=0$, $\omega^3(t)=0$, then the Cauchy 
data is of the form (\ref{t-data})--(\ref{n-data}). 
 
Conversely, if the lateral Cauchy data is of the form 
(\ref{t-data})--(\ref{n-data}) 
for $0 \leq t \leq T$, 
and $\omega$ satisfies
the equation (\ref{M-eq}), with vanishing initial 
data, then $\omega(t)$ is a solution to Maxwell's equations, i.e., 
$\omega^0(t)=0$, $\omega^3(t)=0$. 
\end{theorem} 
 
{\em Proof:} 
The first part of the theorem follows 
from the above considerations if
 we show that $\omega(t)$ is sufficiently regular. 

Since $\omega^2 \in C( \R,H^1(\Omega^2M))$ we see that
${\bf n}\omega^2 \in C( \R,H^{1/2}(\Omega^2 \partial M))$ with
$d{\bf n}\omega^2 \in C( \R,H^{-1/2}(\Omega^2 \partial M))$. Furthermore,
as $\delta \omega_t^1(t) = \delta \delta \omega^2(t) =0$, 
\[
{\bf n}\omega_t^1 \in C( \R,H^{-1/2}(\Omega^2 \partial M)), \quad
\delta \omega^2 \in C( \R,H^{-1/2}(\Omega^2 \partial M)),
\]
which verifies (\ref{n-data}). 

To prove (\ref{t-data}) we use Maxwell duality: 
Consider the forms 
 \[ 
 \eta^{3-k}=(-1)^k*\frac 1\alpha\omega^k. 
\] 
Then $\eta=(\eta^0,\eta^1,\eta^2,\eta^3)$
satisfies Maxwell's equations $\eta_t+\tilde \M \eta=0$ where
$\tilde \M$ is the Maxwell operator with metric $g$
and scalar impedance $\alpha^{-1}$. In sequel, we call
Maxwell's equation with these parameters the adjoint
Maxwell equations and the forms $\eta^j$ the adjoint
solution. Now the formula (\ref{n-data}) for adjoint
solution implies  (\ref{t-data}) for $\omega$.

To prove the converse, it suffices to show that 
$\omega^0(t)=0$. Indeed, the claim $\omega^3(t)=0$ 
follows then by Maxwell duality 
described earlier.
From the equations 
\begin{eqnarray} 
 \omega^0_t -\delta\omega^1 &=& 0,\\ 
 \omega^1_t +d\omega^0 -\delta\omega^2, 
&=&0\label{aux} 
\end{eqnarray} 
it follows that $\omega^0$ satisfies the wave equation 
\[ 
 \omega_{tt}^0 +\delta d\omega^0=0. 
\]  
It also satisfies the initial condition  
$\omega^0(0)=0$ and $\omega^0_t(0)=0$
 and, from (\ref{t-data}),  boundary condition ${\bf t}\omega^0 =0$. Thus, 
$\omega^0 = 0$ for $0\leq t \leq T$. 
 
\hfill$\Box$ 
 
The following definition fixes the solution of 
the forward problem considered in this work. 
 
\begin{definition} 
Let $f=(f^0,f^1,f^2)\in C^{\infty}([0,T]; {\bf \Omega}(\partial M))$  
be a smooth boundary source of the form 
(\ref{t-data}), i.e., $f^0=0$, $f^2_t = -df^1$. 
Further, let $R$ be any right inverse of the 
mapping $\bt$. The solution of the initial-boundary 
value problem 
\[ 
\omega_t+ {\cal M}\omega=0,\quad t>0, 
\] 
\[ 
 \omega(0)=\omega_0\in {\bf L}^2(M),\quad \bt\omega = f, 
\] 
is given by 
\[ 
 \omega = Rf + {\mathcal U}(t)\omega_0 
  - \int_0^t{\mathcal U}(t-s)({\mathcal M}Rf(s) 
+Rf_s(s))ds. 
\] 
\end{definition} 
 
We remark that the boundary data $f$ could be 
chosen from a wider class 
$f\in H^{1/2}(\partial M\times [0,T])$. 
 
Theorem \ref{equivalence} motivates the following definition. 
 
\begin{definition} 
\label{27.11.d}
For solution $\omega$ of Maxwell equations
(\ref{apu1})--(\ref{apu2}) we use the
following notations:
\begin{enumerate} 
\item[i.] The lateral Cauchy data for a solution 
$\omega$ of Maxwell's equations with vanishing 
initial data in the interval $0\leq t\leq T$ 
is given by the pair 
\[ 
(\bt\omega^1(x,t),\bn \omega^2(x,t)),\quad 
(x,t)\in\partial M\times [0,T]. 
\] 
\item[ii.] When $\omega$ satisfies initial condition $\omega(0)=0$ the mapping 
\[ 
{\mathcal Z}^T: C^{\infty}_{00}([0,T],\Omega^1(\partial M))\to
 C^{\infty}_{00}([0,T],\Omega^1(\partial M)),
\]
\[
{\mathcal Z}^T(\bt \omega^1)=\bn \omega^2|_{\partial M \times [0,T]},
\] 
is well defined. We call this map the {\em admittance map}. 
\end{enumerate}
\end{definition} 
Here $ C^{\infty}_{00}([0,T], {\it B})$ consists of $C^{\infty}$ functions of $t$ with values in
a space ${\it B}$, i.e. ${\it B} =\Omega^1(\partial M)$ in definition \ref{27.11.d},
which vanish near $t=0$. 
 
 Note that 
in the  
classical terminology for the electric and magnetic fields, 
${\mathcal Z}^T$
 maps the tangential electric field 
$n\times E|_{\partial M\times [0,T]}$
to the tangential magnetic field 
$n\times H|_{\partial M\times [0,T]}$.

The boundary data and the energy of the field inside $M$
are closely related. 
The following 
result, crucial  from the point of view of boundary control, is 
a 
version of the Blagovestchenskii
formula (see \cite{BeBl} for the case of the scalar wave equation).  
Observe that the following theorem is formulated 
for any solutions of the complete system, not only 
for those that correspond to Maxwell's equations. 
 
\begin{theorem}\label{blacho} 
Let $\omega$ and $\eta$ be smooth 
solutions of the complete system (\ref{complete}). Then the knowledge 
of the lateral Cauchy data 
\[ 
 (\bt\omega,\bn\omega),\quad 
(\bt\eta,\bn\eta), \quad 0\leq t\leq 2T,
\] 
is sufficient for the determination of the inner 
products 
\[ 
 (\omega^j(t),\eta^j(s))_{L^2},\quad 0\leq j\leq 3,\; 
 0\leq s,t\leq T 
\] 
over the manifold $M$. 
\end{theorem} 
 
{\em Proof:} The proof is based on the observation 
that, having the lateral Cauchy data of a solution $\omega$, 
we also have access to the forms  
$d\bt \omega$ and $d\bn \omega$ at the boundary.  
On the other hand, 
$\bt$ commutes with $d$ so that 
\[ 
 \bt d\omega^j = d\bt\omega^j,\quad 
 \bn\delta\omega^j = \bt **d\frac 1\alpha*\omega^j = 
 d\bt\frac 1\alpha*\omega^j = d\bn \omega^j. 
\] 
Let us define the function 
\[ 
 F^j(s,t) = (\omega^j(s),\eta^j(t)). 
\] 
From the complete system, it follows that 
\begin{eqnarray}\label{A-fomrmu} 
 (\partial_s^2 -\partial_t^2)F^j(s,t) 
 &=&  (\omega_{ss}^j(s),\eta^j(t))_{L^2}-
  (\omega^j(s),\eta_{tt}^j(t))_{L^2}\\ 
\noalign{\vskip6pt} \nonumber
&=& 
-((d\delta+\delta d)\omega^j(s),\eta^j(t))_{L^2}+
 (\omega^j(s),(d\delta+\delta d)\eta^j(t))_{L^2} \\ 
\noalign{\vskip6pt} \nonumber
 &=& f^j(s,t). 
\end{eqnarray} 
By applying Stokes theorem we obtain further 
that 
\[ 
 f^j(s,t)=\langle \bn\omega^j,\bt  
\delta\eta^j\rangle 
+\langle \bt\omega^j,\bn d\eta^j\rangle 
-\langle \bt \delta\omega^j,\bn\eta^j\rangle 
-\langle \bn d\omega^j,\bt\eta^j\rangle, 
\] 
where we suppressed for brevity the dependence
 of the boundary values on $s$ and $t$. Now the 
complete system implies that 
\[ 
 d\omega^j=-\omega_s^{j+1}+\delta\omega^{j+2}, 
\quad 
 \delta\omega^j = \omega_s^{j-1} + d\omega^{j-2}, 
\] 
and, similarly, 
\[ 
 d\eta^j=-\eta_t^{j+1}+\delta\eta^{j+2}, 
\quad 
 \delta\eta^j = \eta_t^{j-1} + d\eta^{j-2}. 
\] 
A substitution to the above formulas then gives 
\begin{eqnarray*} 
f^j(s,t) &=& \langle\bn\omega^j,\bt\eta_t^{j-1} 
+d\bt\eta^{j-2} \rangle + 
\langle\bt\omega^j,-\bn\eta_t^{j+1} 
+d\bn\eta^{j+2} \rangle \\ 
\noalign{\vskip6pt} 
& & -\langle\bt\omega_s^{j-1}+d\bt\omega^{j-2},\bn\eta^j 
 \rangle -\langle-\bn\omega_s^{j+1}+d\bn\omega^{j+2},\bt\eta^j 
 \rangle, 
\end{eqnarray*} 
where $d$ stands for the exterior derivative on $\p M$.
Hence, $f^j$ is completely determined by the lateral Cauchy 
data. What is more, we have 
\beq\label{B-fo} 
 F^j(0,t)=F^j(s,0)=0,\quad F_s^j(0,t)=F_t^j(s,0)=0. 
\eeq
Hence, we can solve $F(s,t)$ using (\ref{A-fomrmu}) and (\ref{B-fo})
as claimed. 
\hfill$\Box$ 
  
{\bf Remark 3.} If $\omega$ and $\eta$ are solutions to 
Maxwell's equations (\ref{apu1}-\ref{apu2}), the formulas above simplify. 
We have 
\[ 
 f^0(s,t) = f^3(s,t) = 0, 
\] 
and 
\begin{eqnarray*} 
 f^1(s,t) = \langle \bn\omega^2_s,\bt\eta^1\rangle 
- \langle \bt\omega^1,\bn\eta^2_t\rangle,\quad 
 f^2(s,t) = \langle \bn\omega^2,\bt\eta^1_t\rangle 
- \langle \bt\omega^1_s,\bn\eta^2\rangle. 
\end{eqnarray*} 
Then, for $j=1$ the inner product $
(\omega^1(t),\omega^1(t))_{L^2}$
 defines the energy of the electric field. 
Similarly, for $j=2$ the inner product $
(\omega^2(t),\omega^2(t))_{L^2}$ defines the energy of the magnetic field. 

\section{Inverse problem} 
 
The main objective of this chapter is to prove the following 
uniqueness result for the inverse 
boundary value problem. 
 
\begin{theorem}\label{ip} 
Given $\partial M$ and the admittance map 
${\mathcal Z}_T$, $T>8 \diam(M)$, 
for Maxwell's equations, (\ref{apu1})-- (\ref{apu2}), it
is possible to    uniquely reconstruct the Riemannian manifold, 
$(M,g)$ and the scalar wave impedance, $\alpha$. 
\end{theorem} 
 
Observe that, once we know the travel time metric 
$g$ as well as the wave impedance $\alpha$, 
formula (\ref{travel time metric}) gives the 
metrics $g_\mu$ and $g_\epsilon$, 
which correspond to the 
material parameters $\mu$ and $\epsilon$.

The proof of the above result is divided in 
several parts. The first step, 
 which is discussed in the 
next sections,
is to prove necessary boundary 
controllability results. These results are used, 
in a similar fashion as in  \cite{Ku1}, \cite{KKL}, 
to reconstruct the manifold and the travel time 
metric. 

\subsection{Unique continuation results} 
 
\begin{figure}[htbp]
\begin{center}
\includegraphics[width=6cm]{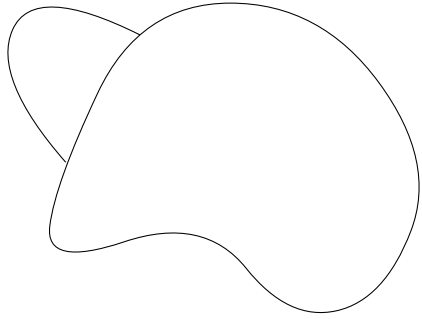} \label{pic 1}
\end{center}
\caption{Manifold $\tilde M$ is obtained by gluing an ``ear'' to $M$.} 
\end{figure}

In the following lemma, we consider extensions of 
differential forms
 outside the manifold $M$. Let $\Gamma\subset\partial M$ be open. 
Assume that 
$\tilde M$ is 
an extension of $M$ across $\Gamma$, i.e.
$M\subset\tilde M$, $\Gamma\subset {\rm int}(\tilde M)$ 
and $\partial M\setminus\Gamma\subset\partial 
\tilde M$. Furthermore, we assume that the metric $g$ and
impedance $\alpha$
are extended smoothly  into $\tilde M$
as $\tilde{g}, \, \tilde{\alpha}$.
In this case, we say that the manifold with scalar impedance
$(\tilde M, \tilde{g},\tilde{\alpha})$ is an {\em extension of $(M,g,\alpha)$ across 
$\Gamma$.} (See Figure \ref{pic 1}). 

{\nottopapertext We have the following simple result.} 
 
\begin{lemma} 
Assume that $\tilde M$ is an extension of $M$ across
an open set $\Gamma\subset \partial M$.
Let $\omega^k$ be
a $k$-form  on $M$ and
$\tilde\omega^k$ be its extension 
by zero  to $\tilde M$. 
Then 
\begin{enumerate} 
\item If $\omega^k\in H(d,\Omega^k M)$ and  
$\bt\omega^k\big|_\Gamma=0$, then $\tilde 
\omega^k\in H(d,\Omega^k\tilde M)$. 
\item If $\omega^k\in H(\delta, 
\Omega^k M)$ and  
$\bn\omega^k\big|_\Gamma=0$, then $\tilde 
\omega^k\in H(\delta,\Omega^k\tilde M)$. 
\end{enumerate} 
\end{lemma} 

{\nottopapertext 
{\em Proof:} 
External differential, in terms of distributions, of $\tilde\omega^k$ can be defined by 
\[ 
 (d\tilde\omega^k,\varphi^k)_{L^2} 
 =(\tilde\omega^k,\delta \varphi)_{L^2},
\] 
where $\varphi^{k}\in \Omega^k\tilde M^{{\rm int}}$ is arbitrary.
However,  by the formula (\ref{stokes1}),  
\[ 
 (\tilde\omega^k,\delta\varphi^k)_{L^2(\tilde M)}= 
 (\omega^k,\delta\varphi^k)_{L^2(M)}=  (d\tilde\omega^k,\varphi^k)_{L^2(M)} + 
  \langle\bt\omega^k,\bn\varphi^k\rangle. 
\] 
Moreover, since $ {\rm supp}\,(\varphi^k) \subset \subset  \tilde M^{{\rm int}}$,
then ${\rm supp}(\bn\varphi^k) \subset \Gamma$,
where $\bt\omega^k$ vanishes. Thus,
\[ 
(\tilde\omega^k,\delta\varphi^k)_{L^2(\tilde M)}
= (d\omega^k,\varphi^k)_{L^2},
\]
i.e.
 $d\tilde\omega^k$ is the zero extension of $d\omega^k$.
In particular,  $d\tilde\omega^k \in L^2(\tilde M)$, 
so that 
$\tilde 
\omega^k\in H(d,\Omega^k\tilde M)$.

The claim concerning the codifferential is proved 
by a similar argument.\hfill$\Box$ 
}
 
As a consequence of this result, we obtain the 
following. 
 
\begin{theorem} 
Let $\omega\in C^1(\R,{\bf L}^2)\cap C(\R,{\bf H}), \, $ ${\bf t}\omega|_{\Gamma \times [0,T]}=0,
\,  {\bf n}\omega|_{\Gamma \times [0,T]}=0$,
be a solution of the equation $\omega_t+{\mathcal M}\omega 
=0$ in $M\times [0,T]$. Let $\tilde\omega$ 
be its extension by zero across $\Gamma\subset 
\partial M$. Then 
the extended form, $\tilde{w}(t)$ satisfies the complete Maxwell's system
on $(\tilde M, \tilde g, \tilde{\alpha})$, i.e. 
$\tilde\omega_t+{\widetilde {\mathcal M}}
\tilde\omega=0$ in 
$\tilde M\times [0,T]$.
\end{theorem} 
 
We are particularly interested in the solutions 
of Maxwell's equations. The following 
result is not directly needed but we have included it,
 since the basic idea is useful when we will
prove the main result of this section.

\begin{lemma}
\label{2.12l}
Assume that $\omega$ in the above theorem satisfies 
Maxwell's equations, i.e., $\omega^0=0$  
and $\omega^3=0$, and $\omega(x,0)=0$.  
If $\bt\omega^1=0$ and 
$\bn\omega^2=0$ on $\Gamma\times [0,T]$, 
then $\omega$  
satisfies Maxwell's equations in the extended 
domain $\tilde M\times [0,T]$. 
\end{lemma} 
 
{\em Proof:} From Theorem \ref{equivalence} it follows 
that, since $\omega$ satisfies Maxwell's equations,  
\[ 
 \bt \omega =  
(0,\bt\omega^1,-\int_0^td\bt\omega^1dt') =0,\quad
 \bn \omega  =  
 (0,\bn\omega^2,\int_0^td\bn\omega^2 dt') =0
\] 
in $\Gamma\times[0,T]$. 
Therefore, the previous theorem shows that  
the  continuation by zero across
$\Gamma, \,$  $\tilde \omega (t)$, satisfies the complete system in 
$\tilde M\times [0,T]$.

However, $\tilde \omega^0(t)=0$,  $\tilde \omega^3(t)=0$ in 
$\tilde M\times [0,T],$ i.e., $\tilde \omega(t)$ satisfies 
Maxwell's equations 
with vanishing initial data 
in the extended manifold  $\tilde M$.
\hfill$\Box$ 

When we deal with a general solution to  Maxwell's equations,
(\ref{apu1})--(\ref{apu2}), which may not satisfy zero initial conditions, 
and try to extend them by zero across $\Gamma$, the arguments of 
Lemma  \ref{2.12l} fail.  Indeed, if $\omega_0 \neq 0$, then
(\ref{2.12.1}) show that $\bn \omega^2 =0$ is not sufficient
for $\bn \omega^1 =0$.
 However, by 
differentiating with respect to time, the 
parasite term $\bn\omega^1(0)$ vanishes. 
This is the motivation why,  in the following 
theorem, we consider the time 
derivatives of the weak solutions.

\begin{figure}[htbp]
\begin{center}
\includegraphics[width=4cm]{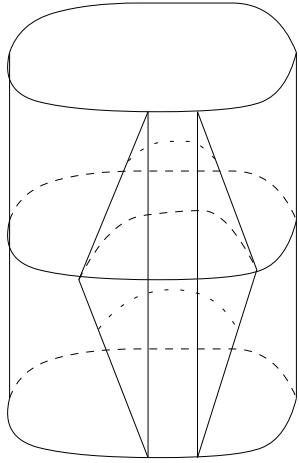} \label{pic 2}
\end{center}
\caption{Double cone of influence.} 
\end{figure}

 
Denote by $\tau(x,y)$  the geodesic distance 
between $x$ and $y$ on $(M,g)$. 
Let $\Gamma\subset\partial M$ be open and $T >0$.
 We use the notation 
\[ 
 K(\Gamma,T)=\{(x,t)\in M\times [0,2T]\mid 
 \tau(x,\Gamma)<T-|T-t|\} 
\] 
for the double cone of influence with base on the slice $t=T$. 
(see Figure \ref{pic 2}.)

\begin{theorem}\label{UCP} 
Let $\omega(t)$ be a weak solution of 
Maxwell's system in the sense of Definition 
\ref{weak solution}  with $\omega_0 =(0,\omega_0^1, 
\omega^2_0,0)$. Assume, in addition, that $\delta\omega^1_0=0$, 
$d\omega_0^2=0$ and $\rho =0$. 
If $\bn\omega^2=0$ in $\Gamma\times]0,2T[$, then 
$\partial_t\omega=0$ in the double cone $K(\Gamma,T)$. 
\end{theorem} 
 
{\em Proof:} 
 Let $\psi \in C^{\infty}_0([-1,1]), \, \int_{-1}^1 \psi(s)ds=1$ be a Friedrich's mollifier.
Then, for any $\sigma >0$ and $\omega(t) \in C(]0,2T[), {\bf L}^2(M))$
satisfying conditions of the Theorem, denote by $\omega_{\sigma}(t)$ its
time-regularization, 
\[
\omega_{\sigma} = \psi_{\sigma} * \omega, \quad 
\psi_{\sigma}(t) =  (1/\sigma) \psi(t/\sigma).
\]
Then $\omega_{\sigma} \in C^{\infty}([\sigma, 2T-\sigma[, {\bf L}^2(M))$ 
continue to be weak solutions to the Maxwell system and, moreover, to Maxwell's equations 
(\ref{apu1})--(\ref{apu2}). Thus, 
\[
{\mathcal M} \omega_{\sigma}  = -\partial _t \omega_{\sigma}  
 \in C^{\infty}([\sigma, 2T-\sigma[, {\bf L}^2(M)),
\]
i.e. $\omega_{\sigma}  \in C^{\infty}([\sigma, 2T-\sigma[, {\mathcal D}({\mathcal M}_{{\rm e}}))$.
Repeating these arguments,
\[
\omega_{\sigma}  \in C^{\infty}([\sigma, 2T-\sigma[, {\mathcal D}({\mathcal M}^{\infty}_{{\rm e}})),
\]
with $ {\mathcal D}({\mathcal M}^{\infty}_{{\rm e}})= \bigcap _{N>1}
 {\mathcal D}({\mathcal M}^N_{{\rm e}})$.

As $\bn \omega_{\sigma}  = \psi_{\sigma}* \bn \omega$,
\[
\bn (\omega_{\sigma})^2  = 0 \quad {\rm on} \,\, \Gamma \times [\sigma, 2T-\sigma[.
\]
Applying  (\ref{2.12.1}),  we see that $\bn \partial_t\omega_{\sigma}  = 0$ on
$\Gamma \times [\sigma, 2T-\sigma[.$

Denote by $\tilde \omega$ the extension by zero of $\omega$ across
$\Gamma$ and 
$\tilde{\eta}_{\sigma}$ that of $\partial_t\omega_{\sigma}$. 
We claim that, in the distributional sense,  
$\tilde{\eta}_{\sigma}$ satisfies the 
complete Maxwell system, for $\sigma<t<2T-\sigma$.
Indeed, let $\varphi=(\varphi^0,\varphi^1, 
\varphi^2,\varphi^3)\in C^\infty_0(]\sigma, 2T-\sigma[, {\bf \Omega}\tilde M^{{\rm int}})$ 
be a test form. 
Using the brackets $[\;\cdot\;,\;\cdot\;]$ to denote 
the distribution duality, that extends the inner 
product 
\[ 
 [\psi,\phi]=\int_0^{2T}(\psi(t),\phi(t))_{{\bf L}^2(\tilde M)}dt, 
\] 
we have 
\begin{eqnarray*} 
[\partial_t\tilde\eta_{\sigma}+{\widetilde {\mathcal M}}\tilde\eta_{\sigma},\varphi] &=& 
 -[\tilde\eta_{\sigma},{\widetilde {\mathcal M}}\varphi +\varphi_t] \\ 
\noalign{\vskip6pt} 
= [\tilde\omega_{\sigma}, 
{\widetilde {\mathcal M}}(\varphi_t) + 
\varphi_{tt}]
 &=&
[\omega_{\sigma}, 
 {\mathcal M}(\varphi_t) + 
\varphi_{tt}]. 
\end{eqnarray*} 
As $\bt \omega_{\sigma} =0$, it follows from the Stokes' theorem and
the fact that $\omega_{\sigma}$ satisfies Maxwell's equations, that
\begin{eqnarray*} 
[\omega, 
{\mathcal M}\varphi_t + 
\varphi_{tt}] = 
 \int_0^{2T}(\omega_{\sigma},{\mathcal M} 
\varphi_t + 
\varphi_{tt})_{{\bf L}^2(M)} dt
= \int_0^{2T}\langle\bn\omega_{\sigma}.
\bt\varphi_t\rangle dt, 
\end{eqnarray*} 
As ${\rm supp}(\bt \varphi) \subset \Gamma \times ]\sigma, 2T-\sigma[$, where
$\bn \omega_{\sigma} =0$, the right side of this equation equals to $0$. In addition,
$\tilde{\bt} \tilde \omega_{\sigma} =0$ for $ t \in ]\sigma, 2T-\sigma[$, where 
$\tilde{\bt}$ is the tangential component on $\partial \tilde M$. Thus, the claim follows.

However, 
$
 \tilde\eta_{\sigma}\in C^\infty(]\sigma,2T-\sigma[,{\bf L}^2(\tilde M)).
$
Therefore, similar considerations to the above shows that this implies that 
\[
\tilde \eta_{\sigma}  \in C^{\infty}([\sigma, 2T-\sigma[, {\mathcal D}^{\infty}({\mathcal M}_{{\rm e}})),
\]
i.e. $\tilde \eta_{\sigma} $ is infinitely smooth in $\tilde M^{{\rm int}} \times [\sigma, 2T-\sigma[$.
Since $\tilde\eta_{\sigma}=0$ outside $M\times \R$, 
the unique continuation result of 
Eller-Isakov-Nakamura-Tataru \cite{EIsNkTa}, that is based on result of
Tataru  \cite{Ta1},\cite{Ta3}

for smooth solutions, implies that $\tilde\eta_{\sigma}=0$ 
in the double cone 
$ 
 \tilde{\tau}(x,\tilde M \setminus M)<T - \sigma-|T-t|, \quad x \in \tilde M,
$
where $ \tilde{\tau}$ is the distance on $(\tilde M, \tilde g)$.
As $\tilde\eta_{\sigma}= \partial_t \omega_{\sigma}$ in $M$, this implies that
$\partial_t \omega_{\sigma}=0$ in the double cone
\beq
\label{3.12.1}
\tau(x,\Gamma)<T - \sigma-|T-t|, \quad x \in  M.
\eeq
When $\sigma \to 0$,  $\tilde\eta_{\sigma}\to\partial_t \omega$, 
in the distributional sense, while the cone (\ref{3.12.1}) tend to
$K(\Gamma,T)$ and the claim of the theorem 
follows.\hfill$\Box$

 We note that the unique continuation result of \cite{EIsNkTa} is 
related to scalar $\epsilon$, $\mu$. However, it is easily generalized 
to the scalar impedance case due to the single velocity of the wave
propagation.

Following the proof of Theorem \ref{UCP}, we can show the 
following variant of Theorem \ref{equivalence}.

\begin{corollary}\label{1.7+2.5} 
Let $\omega(t)$ be a weak solution to the complete Maxwell system in the sense 
of definition \ref{weak solution}, with $\rho=0$,
and, in addition, (\ref{n-data})
 on
$\Gamma \times ]0,T[.$
If $T> 2\,\hbox{diam}(M)$, 
then $\omega^0(t) =0, \, \omega^3(t)=0$ and $\omega(t)$ is a solution of
Maxwell's system for $0<t<T$.
\end{corollary}

{\em Proof:} 
We will consider only $\omega^0$ using the n Maxwell duality for $\omega^3$.
 By remark 1 and (\ref{t-data}),
\beq
\label{wave}
\omega^0_{tt} +\delta d \omega^0=0, \quad {\bf t}\omega^0|_{\partial M \times [0,T]} =0.
\eeq
Also 
\[
\omega^1_t +d\omega^0 -\delta \omega^2 =0,
\]
imply, together with  (\ref{n-data}), that
\[
{\bf n} d \omega^0  = {\bf n} \delta \omega^2 - {\bf n}  \omega^1_t=
d{\bf n}  \omega^2- {\bf n}  \omega^1_t =0
\]
on $\Gamma \times [0,T]$. Together with the boundary condition in (\ref{wave}), 
this shows that the lateral Cauchy data of $\omega(t)$ vanishes on
$\Gamma \times [0,T]$. Using now the wave equation in (\ref{wave}), 
this imply that, due to Tataru's unique continuation 
\cite{Ta1}, \cite{Ta3}, $\omega_0 = 0$ in the double cone $K(\Gamma,T)$.
As $T > 2\hbox{diam}(M)$, this yield that $\omega^0(T/2) = \omega_t^0(T/2)=0$.
It now follows from (\ref{wave}) that
 $\omega^0(t) =0$ for $0<t<T$. 
\hfill$\Box$ 
 
\subsection{Introduction for controllability} 
 
In this section we derive the controllability 
results for the Maxwell system. We divide these results 
in {\em local results}, i.e., controllability of the 
solutions at short times and in {\em global results}, 
where the time of control is long enough so that 
the controlled  electromagnetic
waves fill the whole manifold. Both types of 
results are based 
on the unique continuation of Theorem \ref{UCP} and representation of
inner products of electromagnetic fields over $M$, in a time slice,
 in terms of
integrals of the lateral Cauchy data 
over the boundary $\partial 
M$ over a time interval which is given by Theorem \ref{blacho} . 
 
Consider the initial boundary value problem 
\begin{equation}\label{maxwell eq} 
 \omega_t +{\mathcal M}\omega =0, 
 \quad t>0, 
\end{equation} 
with the initial data $\omega(0)=0$ and 
the electric boundary data of Maxwell type, 
\begin{equation}\label{maxwell data} 
 \bt\omega = (0,f,-\int_0^tdf(t')dt'), 
\end{equation} 
where we assume that $f\in C^\infty_0 
(\R_+, \Omega^1\partial M)$. By Theorem \ref{equivalence}, 
we know that $\omega^0(t)=0$ and $\omega^3(t)=0$. 
 
Let $\tilde\omega$ denote the weak solution of 
Definition \ref{weak solution} with $\rho=0$ and $\tilde \omega(T)=\omega_0$.
Assume, in addition, that 
the conditions of Lemma \ref{weak is maxwell} are 
satisfied so that $\tilde\omega$ satisfies also 
Maxwell's equations. 
 
As we have seen, Stokes formula implies the 
identity 
\begin{equation}\label{control identity} 
(\omega(T),\omega_0) 
=-\int_0^T\langle\bt\omega,\bn\tilde\omega 
\rangle dt. 
\end{equation} 
We refer to this identity as the {\em control identity} 
in the sequel. 
 
\subsection{Local controllability}\label{local controllability  
section} 
 
In this section, we study differential $1-$forms in $M$ 
 that can 
be generated
 by using appropriate boundary sources 
active for short periods of time. 
Instead of a complete characterization of these 
forms, 
we show that 
there is a large enough subspace
in $ L^2(\Omega^1M)$ that can be  
produced by boundary sources.
The difficulty that prevents a complete 
characterization is related to the topology of the 
domain of influence, which can be very complicated. 
 
Let $\Gamma\subset\partial M$ be an open subset of the 
boundary and $T>0$ arbitrary. We define the {\em domain of 
influence} as 
\[ 
 M(\Gamma,T)=\{x\in M\mid \tau(x,\Gamma)< T\}, 
\] 
where $\tau$ is the distance with respect to the 
travel time metric $g$.  
Observe that $M(\Gamma,T) = K(\Gamma,T) \cap \{t=T\}.$  
 
Furthermore, let $\omega$ 
be the strong solution of the initial-boundary 
value problem 
\[ 
 \omega_t+{\mathcal M}\omega =0,\quad \omega(0)=0, 
\] 
with the boundary value 
\[ 
 \bt\omega = (0,f,-\int_0^t df(t')dt'), 
\] 
where $f\in C^\infty_0(]0,T[,\Omega^1\Gamma)$
with $C^\infty_0(]0,T[,\Omega^1\Gamma)$ being a subspace of  forms in
$C^\infty_0(]0,T[,\Omega^1\partial M)$  with support in $\Gamma$.
To emphasize the dependence of $\omega$ on $f$, we 
write occasionally  
\[ 
\omega=\omega^f = (0,(\omega^f)^1,(\omega^f)^2,0). 
\] 
We denote 
\beq
\label{3.12.2}
 X(\Gamma,T)={\rm cl}_{L^2}\{ (\omega^f)^1(T)\mid f\in 
 C^\infty_0(]0,T[, \Omega^1\Gamma )\}, 
\eeq 
i.e., $X(\Gamma,T)$ is the $L^2$--closure of 
the set of the electric fields that 
are generated by $C^\infty_0$--boundary sources on 
$\Gamma\times ]0,T[$.  
Furthermore, we use the notation 
\[ 
H(\delta,M(\Gamma,T)) =\{ 
\omega^2\in H( 
\delta,M),\;{\rm supp}\,(\omega^2) 
\in \overline{M(\Gamma,T)}\}. 
\] 
We will prove  the following result. 
 
\begin{theorem}\label{local control th} 
The set $X(\Gamma,T)$ satisfies 
\[ 
 \delta H^1_0(\Omega^2M(\Gamma,T))\subset 
 X(\Gamma,T)\subset{\rm cl}_{L^2}\bigg(\delta H( 
\delta,M(\Gamma,T) 
)\bigg).
\] 
\end{theorem} 

Here $H^1_0(\Omega^2S), \, S \subset M$ is a subspace of $H^1_0(\Omega^2M)$
of forms with support in ${\rm cl}(S)$.
 
{\em Proof:} The right inclusion is straightforward: 
Since $\omega$ satisfies Maxwell's equations, we have 
$\omega^0(t) = 0$ and, hence, 
\[ 
 \omega^1(T) = \int_0^T\delta\omega^2dt 
\in\delta H(\delta,M(\Gamma,T)). 
\] 
 
To prove the left inclusion, 
we show that any field of the form $\nu^1 
=\delta\eta^2$ with $\eta^2 
\in H^1_0(M(\Gamma,T))$ is in $(X(\Gamma,T)^\perp) 
^\perp$. To this end, let us first 
assume that 
$\omega_0^1\in L^2(\Omega^1M)$ is 
a 1--form such that 
\[ 
 (\omega_0^1,\omega^1)_{L^2}=0 
\] 
for all $\omega^1=(\omega^f(T))^1$ 
generated by    boundary sources 
$f\in C^\infty_0(]0,T[,\Omega^1\Gamma)$. 
Since $\omega^1=\delta\omega^2$, 
it suffices to consider only those forms $\omega_0^1$ 
that are of the form $\omega_0^1=\delta 
\eta^2_0$ for some $\eta^2_0\in H(\delta)$.  
Indeed, by Hodge decomposition (see \cite{Sc}) in ${\bf L}^2(M)$, we have
\ba
\omega^1_0=\hat \omega^1_0+\delta \eta^2_0,
\ea
where $d\hat \omega^1_0=0$, $\bt \hat \omega^1_0=0$
so $\hat \omega^1_0\perp \omega^1$ automatically.
 
Let $\tilde\omega$ be a weak solution, at the time interval $[0,T]$, of the 
initial boundary 
value problem (\ref{ibvp}) with 
$\bt\tilde\omega=0$,  and
\[ 
 \tilde\omega(\;\cdot\;,T)= (0,\omega_0^1,0,0)=\omega_0.
\]  
By our assumption,
\[ 
 (\omega(T),\omega_0)_{{\bf L}^2} 
 = (\omega^1(T),\omega_0^1)_{L^2} =0, 
\] 
and thus, by the control identity, 
(\ref{control identity})
and conditions (\ref{t-data}), (\ref{n-data}), 
\[ 
 \int_0^T\langle \bt \omega,\bn\tilde\omega\rangle 
 =\int_0^T\langle \bt\omega^1,\bn\tilde\omega^2 
\rangle = \int_0^T\langle f,\bn\tilde\omega^2 
\rangle=0, 
\] 
for all 
differential 1-forms 
$f \in C^\infty_0(]0,T[, \Omega^1\Gamma )$. 
Thus,  we have  
\[ 
 \bn\tilde\omega^2=0 \mbox{ on $\Gamma 
\times ]0,T[$}. 
\] 
Furthermore, it is easy to see that, for 
$T+t\in [T,2T]$, we have 
\[ 
 \tilde\omega(T+t)=(0,\tilde\omega^1(T-t), 
 -\tilde\omega^2(T-t),0), 
\] 
and, therefore, also 
\[ 
 \bn\tilde\omega^2=0 \mbox{ on $\Gamma 
\times ]T,2T[$}. 
\] 
But this implies that,  as a distribution, $\bn\tilde\omega^2$ 
vanishes on the whole interval $]0,2T[$ 
 since it is in  
$L^2_{\rm loc}(\R,H^{-1/2}(\partial M))$. 
By applying the  Theorem  
\ref{UCP}, we can deduce that $\tilde\omega_t=0$ 
in the double cone $K(\Gamma,T)$. In particular, 
we have that $d\omega^1_0 = \tilde{\omega}^2_t(T)=0$ in 
$M(\Gamma,T)$. 
 
Let now $\nu^1=\delta\eta^2 
\in\delta H^1_0(\Omega^2M(\Gamma,T))$. Then 
\[ 
 (\nu^1,\omega^1_0)_{L^2} 
 = (\eta^2,d\omega^1_0)_{L^2}=0. 
\] 
This holds for arbitrary $\omega_0^1\in  
X(\Gamma,T)^\perp$, i.e., $\nu \in 
(X(\Gamma,T)^\perp)^\perp =X(\Gamma,T)$. 
\hfill$\Box$ 
 
{\bf Remark 4.} Later in this work, we are mainly interested 
in controlling the time derivatives of 
 electromagnetic  fields.
Let 
us denote 
\[ 
 \Cnull(\Gamma, T)= \{\int_0^t f(t')dt'\mid f\in C^\infty_0(]0,T[,
\Omega^1\Gamma) \}. 
\] 
With this notation,  we have 
\[ 
 X(\Gamma, T) = {\rm cl}_{L^2}\{(\omega^f_t(T))^1\mid f\in\Cnull(\Gamma,T)\}. 
\] 
Indeed, if $\omega^1 =(\omega^f)^1\in X(\Gamma,T)$, then 
$(\omega^f)^1 = (\omega^F_t)^1$, where 
\[ 
 F(t)=\int_0^t f(t')dt'. 
\] 
Conversely, the time derivative of a field $\omega^f$, $f\in 
\Cnull(\Gamma,T)$ satisfies the initial-boundary value problem with 
the boundary source $f_t\in C_0^\infty(]0,T[, \Omega^1\Gamma)$. 

 \subsection{Global controllability}  
 
We start by introducing some notations. Let $\omega$ 
be the strong solution of the initial-boundary 
value problem 
\[ 
 \omega_t+{\mathcal M}\omega =0,\quad \omega(0)=0, 
\] 
with the boundary value 
\[ 
 \bt\omega = (0,f,-\int_0^t df(t')dt'), 
\] 
where $f\in C^\infty_0(]0,T_0[, \Omega^1\Gamma)$, 
$T_0>0$ and 
$\Gamma\subset\partial M$ is an open subset.  
For $T\geq T_0$, we define 
\beq
\label{3.12.3}
 Y(\Gamma,T) = \{\omega^f_t(T)\mid 
 f\in C^\infty_0(]0,T_0[, \Omega^1\Gamma)\}. 
\eeq 
 For $\Gamma=\partial M$ we 
denote $Y(T)=Y(\partial M,T)$.
Our objective  
is to give a characterization of the 
set $Y(T)$ for 
$T_0$ large enough. 
 In the following,
\beq
\label{3.12.4}
{\rm rad}(M)=\max_{x\in M} \tau(x,\p M).
\eeq
 We prove the following 
result. 
 
\begin{theorem}\label{global control th} 
Assume that $T_0> 2{\rm rad}(M). $
Then, 
for $T\geq T_0$, ${\rm cl}_{{\bf L}^2(M)}Y(T)$ is independent of $T$, 
i.e. ${\rm cl}_{{\bf L}^2(M)}Y(T)=Y$, and, moreover, 
\beq
\label{3.12.5}
Y = \{0\}\times\delta H(\delta)\times 
 d\Hnull(d)\times\{0\}. 
\eeq
\end{theorem} 
 
\noindent
{\bf Remark 5.} The result holds also for $Y$ replaced with 
$Y(\Gamma,T)$, when 
\ba
T_0>2\max_{x\in M} \tau(x,\Gamma).
\ea

{\em Proof:} Let $\omega=\omega^f$ be a solution. 
As $f=0$ for $T\geq T_0$, we have 
\[ 
 \bt\omega^1(T)=0, 
\] 
and, consequently,
\begin{eqnarray*} 
 \omega_t(T) &=& -{\mathcal M}\omega(T) 
 = (0,\delta\omega^2(T),-d\omega^1(T),0) \\ 
\noalign{\vskip4pt} 
&\in&\{0\}\times\delta H(\delta)\times 
 d\Hnull (d)\times\{0\}. 
\end{eqnarray*} 
To prove the converse inclusion, we show that the 
space $Y(T)$ is dense in $\{0\}\times\delta H(\delta)\times 
 d\Hnull (d)\times\{0\}$. To this end, 
let $\omega_0 \in\{0\}\times\delta H(\delta)\times 
 d\Hnull (d)\times\{0\}$ and 
  $\omega_0\perp Y(T)$. 
This means that,  for  
arbitrary $\omega=\omega^f$ satisfying the  
initial-boundary value problem (\ref{complete}),  
\begin{equation}\label{orthogonality} 
(\omega_0,\omega_t(T))_{{\bf L}^2} =(\omega_0^1,\omega_t^1(T))_{L^2} 
+(\omega_0^2,\omega_t^2(T))_{L^2} =0. 
\end{equation} 

Let $\tilde{\omega}$ denote the weak solution of the problem 
\ba
& & \tilde{\omega}_t + {\mathcal M}\tilde{\omega} =0, \\
& &
 \bt \tilde{\omega} =0,\quad \tilde{\omega}(T)=\omega_0. 
\ea
Observe that the initial value $\omega_0$ satisfies 
\[ 
 \delta\omega_0^1 =0,\quad d\omega_0^2 =0, 
\] 
which implies that
 $\tilde{\omega}$ satisfies Maxwell's equations. 
Consider the function $F:\R\to\R$, 
\[ 
  F(t) = (\tilde{\omega}(t),\omega_t(t))_{{\bf L}^2}. 
\] 
We have, by using Maxwell's equations, that 
\begin{eqnarray*} 
 F_t(t)&=& (\tilde{\omega},\omega_{tt})_{{\bf L}^2} + (\tilde{\omega}_t, 
\omega_t)_{{\bf L}^2} \\ 
\noalign{\vskip4pt} 
 &=& -(\tilde{\omega}^1,\delta d\omega^1)_{L^2} - (\tilde{\omega}^2, 
 d\delta\omega^2)_{L^2} +(d\tilde{\omega}^1,d\omega^1)_{L^2} + (\delta 
\tilde{\omega}^2,\delta\omega^2)_{L^2}, 
\end{eqnarray*} 
and further, by using Stokes' theorem, 
\[ 
 F_t(t) = -\langle\bt \tilde{\omega}^1(t),\bn d\omega^1(t)\rangle 
 -\langle\bn \tilde{\omega}^2(t),\bt\delta\omega^2(t)\rangle. 
\] 
However, ${\bt}\tilde \omega =0$ and $\delta \omega^2 =\omega^1_t$. Thus,
\[  
 F_t(t) = -\langle\bn \tilde{\omega}^2,\bt\omega^1_t\rangle 
 =-\langle\bn \tilde{\omega}^2,f_t\rangle. 
\] 
On the other hand, the initial condition $\omega(0)=0$, 
together with the orthogonality condition  
(\ref{orthogonality}), imply 
that $F(0)=F(T)=0$,  so that
\[ 
 \int_0^T\langle\bn \tilde{\omega}^2,f_t\rangle dt = -\int_0^T 
 F_t(t)dt =0. 
\] 
Since 
 $\, f \in C^{\infty}_0(]0,T[, \Omega^1 \Gamma)$ is arbitrary, this implies
 that 
\[ 
 \bn \tilde{\omega}^2_t =0 \mbox{ in $\Gamma\times]0,T[$}. 
\] 
But now Theorem \ref{UCP} implies that $\tilde{\omega}_{tt}$ 
vanishes in the double cone $K(\Gamma, T/2)$. 
By the assumption $T_0> 2{\rm rad}(M)$, 
this double cone contains a cylinder 
\[ 
 C = M\times]T/2-s,T/2+s[ 
\] 
\begin{figure}[htbp]
\begin{center}
\psfrag{1}{$\Gamma$}
\includegraphics[width=4cm]{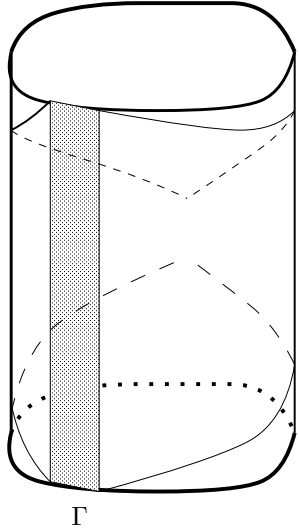} \label{pic 3}
\end{center}
\caption{The double cone contains a slice $\{T/2\}\times M$ and
the waves vanish near the slice $t=T/2$.}
\end{figure}
with some $s>0$ (See Figure \ref{pic 3}). Therefore, 
$\tilde{\omega}_{tt}$ that satisfies Maxwell's equations and 
 homogeneous boundary condition $\bt \tilde{\omega}_{tt}=0$. 
Therefore, it vanishes in the whole $M\times\R$. In particular,  
this means that, with some time-independent 
forms 
$\omega_1$ and $\omega_2$, 
\[ 
 \tilde{\omega}(t)=\omega_1 + t\omega_2, 
\] 
with $\bt\omega_1=0$, $\bt\omega_2=0$.
Again, by Maxwell's equations, we have  
\[ 
 \omega_2 = \omega_t ={\mathcal M}\omega_1 +  
 t{\mathcal M}\omega_2, 
\] 
for all $t$.
Therefore, 
\[ 
 \omega_2 ={\mathcal M}\omega_1,\quad 
{\mathcal M}\omega_2 =0. 
\] 
But then, Stokes' theorem implies that 
\[ 
 (\omega_2,\omega_2)_{{\bf L}^2} = (\omega_2,{\mathcal M}\omega_1)_{{\bf L}^2} 
 = -({\mathcal M}\omega_2,\omega_1)_{{\bf L}^2}=0, 
\] 
i.e., $\omega_2=0$ and  ${\mathcal M}\omega_1=0$. 
Observe that,  by the assumption of the Theorem,
\[ 
 \omega_1 = \tilde{\omega}(T)=\omega_0 = (0,-\delta\nu^2, 
d\nu^1,0)= {\mathcal M}\nu, 
\] 
for some $\nu\in \{0\}\times \Hnull (d)\times H(\delta) 
\times\{0\}$. Therefore, a further application of 
Stokes theorem gives 
\[ 
(\omega_1,\omega_1)_{{\bf L}^2} = (\omega_1,{\mathcal M}\nu)_{{\bf L}^2}
 = -({\mathcal M}\omega_1,\nu)_{{\bf L}^2}=0, 
\] 
i.e., also $\omega_1=\omega_0=0$. The proof is  
therefore complete.\hfill$\Box$ 
 
\subsection{Generalized sources}\label{generalized sources section} 
 
So far, we have treated only smooth boundary 
sources and the corresponding fields. For later 
use, we need more general boundary sources. 
 
Let $Y={\rm cl}_{{\bf L}^2(M)}Y(\partial M,T)$ be the space of  
the time derivatives of electromagnetic fields satisfying Maxwell's 
equations, see (\ref{3.12.3}), (\ref{3.12.5}). We define the wave operator 
\[ 
 W^T:C^\infty_0(]0,T[, \Omega^1\partial M )\to Y, 
\quad f\mapsto \omega_t^f(T), 
\] 
where $T\geq T_0$ and $T_0>2\,{\rm rad}\,(M)$. By means of the wave operator, 
we 
define the ${\mathcal F}-$norm on the space of boundary sources as 
\begin{equation}\label{f norm} 
 \|f\|_{{\mathcal F}} = \|W^T f\|_{{\bf L}^2}. 
\end{equation} 
The definition of this norm is independent of the 
choice of $T\geq T_0$ by conservation of energy. 
 
Notice that by Theorem \ref{blacho}, the knowledge 
of the admittance map ${\mathcal Z}^{2T}$ enables 
us to calculate explicitly the ${\mathcal F}$--norm 
of 
any smooth boundary source. 
 
 To complete the space of boundary sources, 
let us define 
the equivalence $\sim$ of sources by setting 
\[ 
 f\sim g\mbox{ iff $W^Tf=W^Tg$}. 
\] 
Further, we define the space ${\mathcal F}([0,T_0])$ as 
\[ 
 {\mathcal F}([0,T_0])=
C^\infty_0(]0,T[, \Omega^1\partial M )/\sim. 
\] 
Finally, we complete ${\mathcal F}([0,T_0])$  
with respect 
to the norm (\ref{f norm}). Hence, this space, 
denoted by ${\overline{\mathcal F}}([0,T_0])$  
consists of 
Cauchy sequences with respect to the norm 
(\ref{f norm}), denoted as 
\[ 
 \hat f = (f_j)_{j=0}^\infty,\quad f_j\in 
C^\infty_0(]0,T[, \Omega^1\partial M ). 
\] 
Note 
that, for any
$\hat f\in \overline \F$, we can find $\hat h\in \overline \F$
such that
 $\hat h=\hat f$ and $\hat h=(h_j)_{j=1}^\infty$,
 $h_j\in 
C^\infty_0(]\e,T[, \Omega^1\partial M )$ for some $\e>0$.
The reason for this is that Theorem \ref{global control th} 
is valid also with $T_0$ replaced with $T_0-\e$, when $\e$ is 
small enough. Thus, for small $\e>0$, we can define,
for any $\hat f= (f_j)_{j=0}^\infty\in \overline \F$, the translation
 $\hat f(\cdotp+\e)=(f_j(\cdotp+\e))_{j=0}^\infty\in \overline \F$.

These sources are called {\em generalized sources} 
in the sequel. The corresponding electromagnetic waves are denoted 
as 
\beq
\label{gener}
 \omega_t^{\hat f}(t) = \lim_{j\to \infty}\omega_t^{f_j}(t)\quad
\hbox{for }t\geq T_0. 
\eeq
By the isometry of the wave operator, the above limit 
exists in ${\bf L}^2$ for all generalized sources. 

We note that the above construction of the space of generalized
sources in well-known in PDE-control, e.g. \cite{Ru}, \cite{LTr}.

{\bf Remark 6.} Observe that since the wave operator 
$W^T$ is an isometry and  
${\overline{\mathcal F}}([0,T_0])$ was defined 
by closing $C^\infty_0(]0,T[, \Omega^1\partial M )$ 
with respect to the norm (\ref{f norm}), the wave 
operator extends to a one-to-one isometry 
\[ 
 \hat f\mapsto \omega_t^{\hat f}(T),\quad 
{\overline {\mathcal F}}([0,T_0])\to  
{\rm cl}_{{\bf L}^2}(Y(\p M,T)), 
\] 
where the target 
space is completely characterized in the 
previous section.

We say that $\hh\in \overline \F$ is a generalized time derivative of
$\hf\in \overline \F$, if for $T=T_0$,
\beq\label{gen derivative}
& &\lim_{\sigma\to 0+} 
||\frac{\hf(\cdotp+\sigma)-\hf(\cdotp)}\sigma-\hh||_{\overline \F}=
\\ \nonumber
&=&\lim_{\sigma\to 0+} 
||\frac{\omega^\hf_t(T+\sigma)-\omega^\hf_t(T)}\sigma-\omega^\hh_t(T)||_
{{\bf L}^2(M)}=0
\eeq
In this case we denote $\hh={\mathbb D} \hf$, or just $\hh=\p_t\hf$.
In the following, we use spaces 
$
\F^s={\mathcal D}({\mathbb D}^s),\ s \in \Bbb{Z}_+$,
which are spaces of generalized sources that have $s$ generalized
derivatives. Note that, if (\ref{gen derivative}) is valid
for $T=T_0$, it is valid for all $T \geq T_0$ due to the conservation
of $L^2$-norm for Maxwell's equation (energy conservation). Thus, if 
$\hf\in \F^s$, we have, for $T\geq T_0$,
\beq\label{note 2}
\M^s \omega^\hf_t(T)=\p_t^s \omega^\hf_t(t)|_{t=T}\in {\bf L}^2(M).
\eeq
Note 
that $\M$  here is the differential expression given by 
(\ref{M}), rather than an operator with some  boundary conditions. 
Since $\bt \omega^\hf_t=0$ on $\p M\times ]T_0,\infty[$,
we see that  $\bt (\p_t^j\omega^\hf_t)(t)=0$ for $t>T_0$ and $j\leq s-1$.
Thus, for $\hf\in \F^s$ and $T\geq T_0$, 
\beq\label{smoothness of gen. wave}
\omega^\hf_t\in \bigcap_{j=0}^s (C^{s-j}([T,\infty[,{\mathcal D}(\M_e^j))
\cap \hbox{Ran}\,(\M_e)).
\eeq
Moreover, by (\ref{smoothness of gen. wave}) and Lemma 
\ref{lem 2.4},
\ba
\omega^\hf_t(T)\in {\bf H}^s_{loc}(M^{int})\quad\hbox{for }T\geq T_0.
\ea

Next we consider dual spaces to the domains of powers of $\M_e$. Since
${\bf H}_0^s\subset {\mathcal D}({\mathcal M}_{\rm e}^s)$, 
we have $({\mathcal D}({\mathcal M}_{\rm e}^s))'\subset 
{\bf H}^{-s}$. Similarly, we see that
 $ 
{\bf H}^{-s}_0\subset ({\mathcal D}({\mathcal M}_{\rm e}^s))'$.
These facts will be needed later in the construction of focusing sources.

\subsection{Reconstruction of the manifold} 
 
In this section we will show how to
 determine the manifold, $M$ and the travel 
time metric, $g$ from the boundary measurements 
of the admittance map ${\mathcal Z}$.
We will show that the boundary data determines the 
set of {\em boundary distance functions}. 
The basic idea is to use a slicing principle, when 
we control the supports of the waves generated by  
boundary sources.  
 
We start by fixing certain notations.  
Let  times $T_0<T_1<T_2$ 
satisfy
\[ 
 T_0>2\,{\rm rad}(M),\quad T_1\geq T_0+{\rm diam}(M), 
\quad T_2\geq 2 \,T_1
\] 
We assume in this section that the admittance map 
${\mathcal Z}^{T_2}$ is known.

\begin{figure}[htbp]
\begin{center}
\psfrag{1}{$M$}
\psfrag{2}{$t=0$}
\psfrag{3}{$t=T_0$}
\psfrag{4}{$t=T_1$}
\psfrag{5}{$t=T_2$}
\includegraphics[width=4cm]{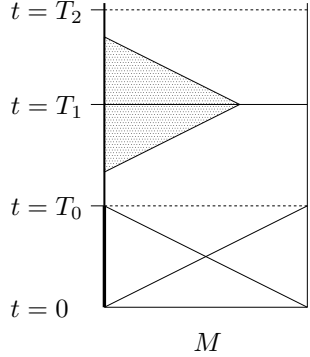} \label{pic 4}
\end{center}
\caption{The sources $\hat f$ of the waves $\omega^\hf(x,t)$ are supported
on the time-interval $[0,T_0]$ which enables us to control the waves
at times $t>T_0$. In the construction of the manifold, the supports of the  waves are
considered at time $t=T_1$. To this end, we use the  unique continuation in
double-cones (triangle in the figure) 
intersecting the boundary in the layer $\p M\times [T_0,T_1]$.
Note that in this layer it is crucial that $\hat f=0$.
}
\end{figure}

Let $\Gamma_j\subset\partial M$ be open 
disjoint sets, $1\leq j\leq J$ and $\tau_j^-$ and 
$\tau_j^+$ be positive times with 
\[ 
 0<\tau_j^-<\tau_j^+\leq {\rm diam}(M), 
 \quad 1\leq j\leq J. 
\] 
We define the set $S = S(\{\Gamma_j,\tau_j^-, 
\tau_j^+\})\subset M$  as an intersection 
of slices, 
\begin{equation}\label{def of S} 
  S =\bigcap_{j=1}^J  \left( M(\Gamma_j,\tau_j^+)\setminus 
M(\Gamma_j,\tau_j^-) \right). 
\end{equation} 
Our first goal is to find out, by boundary measurements, 
whether the set $S$ contains an open ball or not. 
To this end, we give the following definition. 
 
\begin{definition}\label{support sources} 
The set $Z = Z(\{\Gamma_j,\tau_j^-, 
\tau_j^+\})_{j=1}^J$ consists of those generalized sources 
$\hat f\in {\mathcal F}^\infty([0,T_0])$ that produce  
waves $\omega_t=\omega^{\hf}_t$ 
 with
\begin{enumerate} 
\item $\omega_t^1(T_1)\in X(\Gamma_j,\tau_j^+)$, for 
all $j$, $1\leq j\leq J$, 
\item $\omega_t^2(T_1) = 0$, 
\item $\omega_{tt}(T_1)=0$ in $M(\Gamma_j,\tau_j^-)$, 
for 
all $j$, $1\leq j\leq J$. 
\end{enumerate} 
\end{definition}

{\bf Remark 7.} 
Observe that, since $\omega_t$ satisfies Maxwell's 
equations, we have, in particular, 
\[ 
 \omega_{tt}^2 = -d\omega_t^1,\quad 
 \omega_{tt}^1 =\delta\omega_t^2. 
\] 
These identities imply that, at $t=T_1$,  
$\omega_{tt}$ is of the form 
\[ 
 \omega_{tt}(T_1)= (0,0,\omega_{tt}^2(T_1),0) 
 =(0,0,d\eta^1,0), 
\] 
for the 1--form $\eta^1=-\omega_t^1$, and 
\[ 
 {\rm supp}(d \eta^1)\subset S. 
\] 
This observation is crucial later when 
we will  discuss  
focusing waves. 
 
The central tool for reconstruction the manifold 
is the following theorem.

\begin{figure}[htbp]
\begin{center}
\psfrag{1}{$A$}
\includegraphics[width=5cm]{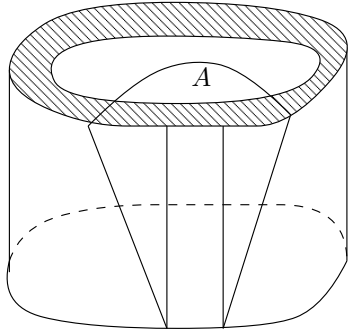} \label{pic 5}
\end{center}
\caption{In Definition \ref{support sources} we can consider e.g.
the case $\Gamma_1=\Gamma$, $\tau^+_1=s_1$,  $\tau^-_1=0$,
and   $\Gamma_2=\p M$, $\tau^+_2=\diam(M)$,  $\tau^-_1=s_2$.
Then the waves that satisfy the definition have the following properties:
By 1., the wave $(\omega^\hf_t)^1(T_1)$ coincides with 
a wave that is supported in $M(\Gamma,s_1)$. This domain of influence
on the figure is the upper part of the cone of influence.
Thus, $d(\omega^\hf_t)^1(T_1)=(\omega^\hf_{tt})^2(T_1)$
 is supported in $M(\Gamma,s_1)$. 
By 2., the wave $(\omega^\hf_{tt})^2(T_1)$ vanish in the boundary layer
$M(\p M,s_2)$. Combining these, we see that  $\omega_{tt}^{\hat f}(T_1)$
is supported in $A=M(\Gamma,s_1)\setminus M(\p M,s_2)$.
}
\end{figure}

\begin{theorem}\label{alternative} 
Let $S$ and $Z$ be defined as above. The following 
alternative holds: 
\begin{enumerate} 
\item If $S$ contains an open ball, then 
${\rm dim}(Z)=\infty$, 
\item If $S$ does not hold an open ball, then 
$Z = \{0\}$. 
\end{enumerate} 
\end{theorem} 
 
In order to prove the above alternative, we need the 
following observability result that will be also useful 
 later. 
 
\begin{theorem}\label{is in Z th} 
Given the boundary map ${\mathcal Z}^{T_2}$, we can 
determine whether a given boundary source $\hat f 
\in{\mathcal F}^\infty([0,T_0])$ is in 
the set $Z$ or not. 
\end{theorem} 
 
{\em Proof:} Let 
$\hat f=(f_k)_{k=0}^\infty \in {\mathcal F}^\infty([0,T_0])$ be a generalized source. 
Consider first the question whether $(\omega_t^{\hat f})^1(T_1) 
\in X(\Gamma_j,\tau_j^+)$. By  Remark 4,  
this is equivalent to the existence of a 
sequence, 
\[
 \hat h=(h_\ell)_{\ell=0}^\infty,\quad h_\ell \in\Cnull(\Gamma_j, 
\tau_j^+), 
\]
such that 
\beq
\label{20.11.1}
 \lim_{k,\ell \to \infty} \| (\omega_t^{f_k})^1(T_1) 
-(\omega_t^{h_\ell})^1(\tau_j^+)\| =0. 
\eeq 
By linearity of the initial-boundary value problem, 
we have 
\[ 
\| (\omega_t^{f_k})^1(T_1) 
-(\omega_t^{h_\ell})^1(\tau_j^+)\| =\|(\omega^{g_{k,\ell}})^1(T_1)\|, 
\] 
where the source $g_{k,\ell}$ is  
\[ 
 g_{k,\ell}(t) = (f_k)_t(t) - (h_\ell)_t(t+\tau_j^+ - T_1)\in C_0^\infty 
( ]0,T_1[;\Omega^1 \p M). 
\] 
However, 
by Lemma \ref{blacho}, $\|(\omega^{g_{k,\ell}})^1(T_1)\|$ is completely
 determined 
by the admittance 
map, ${\mathcal Z}^{T_2}$ making possible to verify (\ref{20.11.1}). 
 
In a similar fashion, Condition 2 of the definition of  $Z$ 
is valid for $\hat f$, if 
\[ 
 \lim_{k\to \infty}\|(\omega_t^{f_k}(T_1))^2\|= 0, 
\]  
and this condition can also be verified via the admittance map, ${\mathcal Z}^{T_2}$.
 
Finally, consider Condition 3. We assume here that we 
already know that $\hat f$ satisfies Conditions 1--2. First, 
we observe that $\omega_{tt}=\omega_{tt}^{\hat f}$ satisfies 
\[ 
 (\partial_t+{\mathcal M})\omega_{tt}=0 \mbox{ in $M\times\R_+$}, 
\] 
along with the boundary condition 
\[ 
 \bt\omega_{tt} =0 \mbox{ in $\partial M\times [T_0,\infty[$}. 
\] 
If Condition 3 holds, by the finite propagation speed, 
$\omega_{tt}$ vanishes in a double cone around $\Gamma_j$, 
i.e., 
\[ 
 \omega_{tt}=0\mbox{ in $K_j=\{(x,t)\in M\times\R_+\mid \tau(x,\Gamma_j) 
+|t-T_1|<\tau_j^-\}$}, 
\] 
for all $j=1,\ldots,J$. In particular, this means that, in each $K_j$, 
$\omega_t$ does not depend on time, and Condition 2 implies that 
$\omega_t^2=0$ in $K_j$. Hence, we have 
\begin{equation}\label{cond 3} 
 \bn\omega_t^2= {\mathcal Z^{T_2}f}=0\mbox{ on $\Gamma_j\times 
 ]T_1-\tau_j^-,T_1+\tau_j^-[$}. 
\end{equation} 
Conversely, assume that condition (\ref{cond 3}) holds together 
with Conditions 1--2. Then $\omega_t$ satisfies 
\[ 
 (\partial_t+{\mathcal M})\omega_t =0\mbox{ in $M\times\R_+$} 
\] 
with the boundary conditions 
\[ 
 \bt\omega_t^1=0,\quad \bn\omega_t^2=0\mbox{ in 
 $\Gamma_j\times]T_1-\tau_j^-,T_1+\tau_j^-[$}. 
\] 
Here we used the fact that $T_1-\tau_j^->T_0$,
so that $\hat{f}=0$ in $\Gamma_j\times]T_1-\tau_j^-,T_1+\tau_j^-[$.   
Now the Unique Continuation Principle,
given by Theorem  \ref{UCP},  
implies that $\omega_{tt}=0$ in $K_j$ and, in 
particular, Condition 3 is valid. The proof 
is complete as it is clear that the condition 
(\ref{cond 3}) is readily observable if the 
admittance map, ${\mathcal Z}^{T_2}, \, T_2 > T_1 + \tau _j,$ 
 is known.
\hfill$\Box$ 
 
Now we can give the proof of Theorem \ref{alternative}. 
 
{\em Proof of Theorem \ref{alternative}:} 
Assume that there is an open ball $B\subset S$. 
Let $0\neq\varphi\in \Omega^2 B$ be an arbitrary smooth 2--form 
with 
 ${\rm supp}\,(\varphi)\subset B$. From the 
global controllability result, Theorem \ref{global 
control th}, it follows the existence of 
a generalized 
source $\hat f \in \overline{{\mathcal F}([0,T_0])}$ 
such that 
\beq
\label{20.11.2} 
 \omega^{\hat f}_t(T_1)=(0,\delta\varphi,0,0). 
\eeq
Moreover, $\varphi \in \Omega^2B$ implies that
$\varphi \in {\mathcal D}({\mathcal M}_e^s)$ for any $s>0$
so that $\hat f \in {\mathcal F}^{\infty}([0,T_0])$. 

We will now show that $\hat f\in Z$. 
Indeed, Conditions 1--2 are obvious from 
the definition (\ref{20.11.2}) of $\hat{f}$. 
Finally, we observe that 
\[ 
 \omega^{\hat f}_{tt}(T_1) = -{\mathcal M} 
\omega_t^{\hat f}(T_1) = (0,0,-d\delta\varphi,0), 
\] 
so Condition 3 is also satisfied. 
This proves the first statement of the theorem. 
 
To prove the second part, assume that $S$ does 
not contain an open ball. Suppose, on the contrary 
to the claim, that there is a non-vanishing 
source $\hat f\in Z$ 
which produces the 
wave $\omega(t) =\omega^{\hat f}(t)$. 
Then, by Conditions 1 and 2 in  Definition 
\ref{support sources},  
\[ 
 {\rm supp}(\omega_t(T_1))\subset 
\bigcap_{j=1}^J M(\Gamma_j,\tau_j^+) = S^+. 
\] 
Furthermore, 
\[ 
 \omega_{tt}(T_1)=-{\mathcal M}\omega_t(T_1), 
\] 
so  that 
 \[ 
 {\rm supp}(\omega_{tt}(T_1))\subset S^+. 
\] 
On the other hand, Condition 3 in Definition 
\ref{support sources} imply that 
\[ 
 \omega_{tt}(T_1)=0\mbox{ in $\bigcup_{j=1}^J 
M(\Gamma_j,\tau_j^-)=S^-$.} 
\] 
Thus ${\rm supp}\omega_{tt}(T_1) \subset S^+\setminus S^-$.
However, 
if the set $S$ does not contain an open ball, then
the set 
$S^+\setminus S^-$ in nowhere dense. 
Since $\omega_{tt}(T_1)$ is smooth in $M^{{\rm int}}$,  it 
 vanishes in $M$. 
In particular, Maxwell's equations 
imply that 
\beq\label{refe +}
 d\omega_t^1(T_1)= -\omega_{tt}^2(T_1) =0. 
\eeq 
On the other hand, $\omega_t\in {\rm cl}_{{\bf L}^2}(Y(\partial M,T_1))$, 
so Theorem \ref{global control th} implies that $\omega_t^1(T_1)$ 
is of the form 
\[ 
 \omega_t^1(T_1) = \delta\eta^2 ,
\] 
for some 2--form $\eta$. 
Since, in addition, $\bt\omega_t^1(T_1)=0$
it then follows, by
 Stokes formula,  that 
\[ 
 (\omega_t^1(T_1),\omega_t^1(T_1))_{L^2} = 
(\delta\eta^2,\omega_t^1(T_1))_{L^2} 
 =(\eta^2,d\omega^1_t(T_1))_{L^2} = 0. 
\] 
Together with Condition 2, this implies 
\[ 
 \omega_t(T_1) =0, 
\] 
contradicting to the assumption
 $\hat f \neq 0$.  
The proof is complete. 
\hfill$\Box$  
 
We are now ready to construct the set of the boundary 
distance functions. For each $x\in M$, 
the corresponding boundary 
distance function, $r_x$ is a continuous function on $\p M$ given by 
\[ 
 r_x: \p M\to\R_+,\quad r_x(z)=\tau(x,z), \quad z \in \partial M. 
\] 
They define {\it the 
 boundary distance map} ${\mathcal R}:M\to C(\p M)$,
${\mathcal R}(x)=r_x$, which is continuous and injective
(see \cite {Ku5}, \cite{KKL}). We shall denote 
the set of all boundary distance functions, i.e., the image
of ${\mathcal R}$, by 
\[ 
 {\mathcal R}(M)=\{r_x\in C(\partial M)\mid x\in M\}. 
\] 
It can be shown (see \cite{Ku5}, {\cite{KKL})
  that, given
 ${\mathcal R}(M) 
\subset L^\infty(\partial M)$ 
we can endow it, in a natural way, 
 with  
a  differentiable structure and a metric tensor $\tilde g$, 
so that $({\mathcal R}(M),\tilde g)$ becomes an isometric
copy of $(M,g)$, 
\[ 
 ({\mathcal R}(M),\tilde g)\cong (M,g). 
\] 
Hence, in order to reconstruct 
the manifold (or more precisely, the isometry type of the manifold),
 it suffices to determine the set,  ${\mathcal R}(M),$ of the boundary 
distance functions. The following result is therefore 
crucial. 
 
\begin{theorem} 
Let the 
admittance map
 ${\mathcal Z}^{T_2}$ be given. Then, for 
any $h\in C(\partial M)$, we can find out whether 
$h\in{\mathcal R}(M)$. 
\end{theorem} 
 
{\em Proof:} The proof is based on  
a discrete approximation process.  
First, we observe that the condition 
$h\in {\mathcal R}(M)$ is equivalent to 
the condition that for any sampling 
$z_1,\ldots,z_J\in\partial M$ of 
the boundary, there must be $x\in M$ 
such that 
\[ 
 h(z_j)= \tau(x,z_j),\quad 1\leq j\leq J. 
\] 
Let us denote $\tau_j = h(z_j)$. By the 
continuity of the distance 
function, $\tau(x,z)$ in $x \in M,\,z \in \p M$,
we 
deduce that the above condition is  
equivalent to the following one: 

\noindent
For any $\varepsilon>0$, the points 
$z_j$ have neighborhoods $\Gamma_j\subset 
\partial M$ with ${\rm diam}(\Gamma_j)<\varepsilon$, 
such that 
\begin{equation}\label{int condition} 
{\rm int}\bigg(\bigcap_{j=1}^J M(\Gamma_j,\tau_j+\varepsilon) 
\setminus M(\Gamma_j,\tau_j-\varepsilon)\bigg)\neq 
\emptyset. 
\end{equation} 
On the other hand, by Theorem \ref{alternative}, 
condition 
(\ref{int condition}) is equivalent to 
\[ 
 {\rm dim}\bigg( Z(\{\Gamma_j,\tau_j+\varepsilon, 
\tau_j-\varepsilon\}\bigg)=\infty, 
\] 
that, by means of  Theorem \ref{is in Z th}, 
 can be verified via boundary data. 
  \hfill$\Box$

As a consequence, we obtain the main result of this section. 
 
\begin{corollary} 
The knowledge of the admittance ${\mathcal Z}^{T_2}$ 
is sufficient for the reconstruction of the manifold, $M$ 
endowed with the travel time metric, $g$.
\end{corollary} 
 
Having the manifold reconstructed, the rest of this article 
is devoted to the reconstruction of the wave impedance, $\alpha$.

\subsection{Focusing sources} 
 
In the previous section it was shown 
 that, using 
boundary data, one can control  supports
 of the  
2--forms $(\omega_{tt}^{\hat{f}})^2(t)$. In this section, the goal is 
to construct a sequence 
of  sources, $(\hat{f}_p), \, p=1,2,\cdots$ such that,
when $p \to \infty$, the corresponding forms $(\omega_{tt}^{\hat{f}_p})^2(T_1)$
 become
 supported at a single 
point,
while $(\omega_{tt}^{\hat{f}_p})^1(T_1)=0$.  For $t \geq T_1$,  
these fields behave like point sources, a fact 
that turns out to be useful for reconstructing 
the wave impedance. 
 
In the following, let $\ud_y$ denote the Dirac delta 
at $y\in M^{{\rm int}}$, i.e., 
\[ 
 \int_M\ud_y(x)\phi(x)dV_g(x)=\phi(y), 
\] 
where $\phi\in C^\infty_0(M)$ and $dV_g$ is volume form of $(M,g)$. 
 
Since the Riemannian manifold $(M,g)$ is already found, 
we can choose $\Gamma_{jp}\subset\partial M$, 
$\,0<\tau_{jp}^-<\tau_{jp}^+<{\rm diam}(M)$, so that 
\beq
\label{20.11.3} 
 S_{p+1}\subset S_p,\quad \bigcap_{p=1}^\infty S_p 
 = \{y\},\; y\in M^{int}. 
\eeq 
Then, 
 $Z_p=Z(\{\Gamma_{jp},\tau_{jp}^-, 
\tau_{jp}^+\}_{j=1}^{J(p)}\})$ is  
the corresponding set 
of generalized sources defined in Definition 
\ref{support sources}.

\begin{definition} 
Let $S_p, \, p=1,2,\cdots,$ satisfy (\ref{20.11.3}).   
 We call the sequence $\tilde f =(\hat f_p), \, p=1,2,\cdots$
 with 
$\hat f_p\in Z_p$, a {\em focusing sequence} 
 of generalized sources of order  $s$  
(for brevity, focusing sources), $ s \in \Bbb{Z}_+$, if 
there is a distribution-form $A=A_y$ on $M$ such that 
\[ 
 \lim_{p\to \infty} (\omega^{\partial_t\hat f_p}_t(T_1) 
,\eta)_{{\bf L}^2}  = (A_y,\eta)_{{\bf L}^2},
\] 
for all $\eta\in {\mathcal D}({\mathcal M}_{\rm e}^s)$. 
\end{definition} 
 
\noindent
{\bf Remark 8.} Observe that, by the identity, 
\begin{equation}\label{omega tt} 
 \omega^{\partial_t \hat f_p}_t = \omega_{tt}^{\hat f_p} 
\end{equation} 
and  Remark 7, 
the 
electromagnetic 
wave $\omega^{\partial_t \hat f_p}_t(T_1)$ is supported in 
${\rm cl}(S_p)$, so $A_y$ must be supported on $\{y\}$. 
 
We will show the following result. 
 
\begin{lemma} 
Let the admittance map ${\mathcal Z}^{T_2}$ be given.
Then, for any $s \in \Bbb{Z}_+$ and any sequence  of generalized sources,
$(\hat f_p), \, p=1,2\cdots$, one can determine if
$(\hat f_p)$ is a focusing sequence or not.  
\end{lemma} 
 
{\em Proof:} Let $\eta\in{\mathcal D}({\mathcal M} 
_{\rm e}^s)$. We decompose  $\eta$ as 
$
 \eta=\eta_1+\eta_2, 
$ 
where 
\[ 
 \eta_1\in {\mathcal D}({\mathcal M} 
_{\rm e}^s)\cap{\rm cl}(Y),\quad \eta_2\in 
{\mathcal D}({\mathcal M} 
_{\rm e}^s)\cap Y^\perp. 
\] 
By the global controllability result, Theorem 
\ref{global control th},  and isometry 
of the wave operator $W^T$, $T\geq T_0$,  
\[ 
 \eta_1 = \omega_t^{\hat h},\quad 
 \hat h\in {\mathcal F}^s([0,T_0]). 
\] 
Since $\omega_t^{\partial_t \hat f_p}\in {\rm cl}_{{\bf L}^2}(Y)$, 
so that $\omega_t^{\partial_t \hat f_p} \perp \eta _2$, 
the condition that $\tilde f$ is a focusing source 
is tantamount to the existence of the limit 
\begin{equation}\label{limit} 
 (A_y,\eta)=\lim_{p\to \infty} (\omega^{\partial_t  
\hat f_p}_t(T_1),\omega_t^{\hat h})_{L^2}, 
\end{equation} 
for all $\hat h\in{\mathcal F}^s([0,T_0])$. However, 
by Theorem \ref{blacho},  the existence of this
limit can be 
verified if we are given ${\mathcal Z}^{T_2}$.
 
Conversely, assume that the limit (\ref{limit}) does 
exist for all $\hat h\in{\mathcal F}^s([0,T_0])$. 
Then, by the Principle of Uniform Boundedness,  
the mappings 
\[ 
 \eta\mapsto(\omega^{\partial_t \hat f_p}_t(T_1) 
,\eta)_{L^2}, \quad p \in \Bbb{Z}_+, 
\] 
form a uniformly bounded family in the dual of 
${\mathcal D}({\mathcal M}_{\rm e}^s)$. By the
Banach-Alaoglu theorem, we find a weak$^*$-convergent 
subsequence 
\[ 
 \omega^{\partial_t \hat f_p}_t(T_1)\to 
 A_y\in \bigg({\mathcal D}({\mathcal M}_{\rm e}^s) 
\bigg)', 
\] 
which is the sought after limit distribution-form.\hfill$\Box$

Since ${\rm supp}(A_y)$ is a point, 
$A_y$ consists of the Dirac delta and its derivatives. The 
role of the smoothness index, $s$ is just to select the 
order of this distribution, as is seen in the 
following result. 
 
\begin{lemma} 
\label{Lm2.15}
Let $ A_y = \lim_{p\to \infty} \omega^{\partial_t  
\hat f_p}_t(T_1)$ is a distribution of order $s=3$.
Then $A_y$ is of the form 
\begin{equation}\label{Ay} 
 A_y(x)= (0,0,d(\lambda\ud_y(x)),0), 
\end{equation} 
where $\lambda$ is a 1--form at $y$, $\, \lambda \in T_y^*M$. 
Furthermore, for any $ \lambda \in T_y^*M$ there is a focusing source
$\tilde{f}=$ with $A_y$ of form (\ref{Ay}). 
\end{lemma} 
 
{\em Proof:}  
From the results in Section \ref{generalized sources section},
 we deduce that, when $s=3$, 
\[ 
 A_y\in\left({\mathcal D}({\mathcal M}_{\rm e}^3) 
\right)'\subset {\bf H}^{-3}. 
\] 
Furthermore, from  Remark 7, the 
electromagnetic waves
 (\ref{omega tt}) are of the form  
\beq\label{refe A}
\omega^{\partial_t\hat f_p}_t(T_1) =(0,0,d\eta_p,0), 
\eeq
for some 1-forms $\eta_p$. 
Combining ( these with the fact that ${\rm supp}(A_y) 
=\{y\}$, we see that 
$
 A_y = (0,0,A^2_y,0).
$ 
Here $A^2_y$, expressed, for example,  in Riemann normal coordinates  
$(x^1,x^2,x^3)$ near  $y$, 
must be of the form 
\[ 
 A^2_y(x)=a^j\ud_y(x) \theta_j + b^{jk} 
\partial_k\ud_y(x) 
 \theta_j, 
\] 
where $\theta_j = (1/2)e_{j k\ell}dx^k\wedge dx^\ell$
 and $e_{j k\ell}$ is the totally antisymmetric permutation symbol. 
Furthermore, by (\ref{refe A}),  
\[ 
 dA^2_y = 
 (a^j\partial_j\ud_y(x) + b^{jk}\partial_k\partial_j 
 \ud_y(x))dV_g=0. 
\] 
Let $\varphi$ be a compactly supported test function 
and, in the vicinity $U$ of $y$, 
\[ 
  \varphi(x) = x^j,\quad j=1,2,3. 
\] 
It follows that 
\[ 
0= (dA^2_y,\varphi)= a^j. 
\] 
Further, let $\psi$ be a compactly supported  
test function 
and, in the vicinity $U$ of $y$, 
\[ 
  \psi(x) = x^jx^k,\quad j,k=1,2,3. 
\] 
As before, we obtain 
\[ 
0=(dA^2_y,\psi)= b^{jk}+b^{kj}. 
\] 
Thus, $b^{jk}$ may be represented as
$b^{jk} =e^{jk\ell}\lambda_\ell, \, \lambda_\ell \in T^*_yM$,
implying that  
\begin{equation}\label{a2} 
 A^2_y(x) = e^{jk\ell}\lambda_\ell\partial_k\ud_y(x) 
\theta_j. 
\end{equation} 
By the properties of the permutation symbols, $e^{jk\ell}$,
\[ 
 e^{jk\ell}\theta_j = \frac 12 e^{jk\ell}e_{jpq}dx^p 
 \wedge dx^q = \delta^k_p\delta^\ell_q dx^p 
 \wedge dx^q. 
\] 
Substituting this expression back to (\ref{a2}),
we finally obtain 
\[ 
 A^2_y(x) = \lambda_\ell\partial_k\ud_y(x) 
 dx^k\wedge dx^\ell 
 = d(\ud_y(x)\lambda_\ell dx^\ell), 
\] 
as claimed.\hfill$\Box$ 
 
By the above 
results, for any $y\in M^{int}$ and $\lambda \in T^*_yM$,
we can, in principle, find 
focusing sequences $\tilde{f}$ such that 
$\omega_{tt}^{\tilde{f}}(T_1) =A_y$, where $A_y$ is of form 
(\ref{Ay}).  We should, however, stress that, at this stage, we 
can not control the corresponding $\lambda=\lambda (y)$.

Consider now a family of focusing sources
 $\tilde{f}_y, \, y \in M^{{\rm int}}$,
with the corresponding $1-$forms $\lambda (y)$. 
 
\begin{lemma} \label{lem: 9.3}
 Given the admittance map ${\mathcal Z}^{T_2}$,  
it is possible to determine 
whether the map $y\mapsto\lambda_y$ is  a  
nowhere vanishing 1-form valued $C^\infty$--function.  
\end{lemma}  
 
{\bf Proof:} 
Let $\varphi \in \Omega^1M^{{\rm int}}$ be an arbitrary 
compactly supported test $1-$form. By Theorem \ref{global control th},   
there is a generalized 
source $\hat h\in{\mathcal F}^\infty$ such that 
\[ 
 (\omega_t^{\hat h})^1(T_1) = \varphi. 
\] 
Let $\tilde f = (\hat f_p), p=1,2,\cdots,$ be a focusing source 
of order $s=3$. Then, by Lemma \ref{Lm2.15} and the 
definition of the focusing sources, we have 
\begin{eqnarray}
\lim_{p\to \infty}(\omega_{tt}^{\hat{ f}_p}(T_1), 
 \omega^{\hat h}(T_1)) &=& (A_y,\omega^{\hat h} 
 (T_1)) \\ 
\noalign{\vskip4pt} 
= (d(\lambda\ud_y),(\omega^{\hat h})^2(T_1)) 
&=& \int_M \lambda_y \ud_y\wedge*\delta(\omega^{\hat h})^2(T_1). 
\nonumber 
\end{eqnarray} 
Further, by Maxwell's equations,  
\begin{eqnarray*} 
 \lambda_y\wedge*\delta(\omega^{\hat h})^2(T_1) 
 &=&\lambda_y\wedge*(\omega^{\hat h}_t)^1(T_1)\\ 
\noalign{\vskip4pt} 
 &=&\lambda_y\wedge *\varphi(y), 
\end{eqnarray*} 
Here, for $\lambda, \, \eta \in T^*_yM$,
\[
\lambda\wedge *\eta = \langle \lambda, \eta \rangle _y \, dV_g =
g^{jk} \lambda _j \eta _k \, dV_g.
\] 
Thus, 
\beq
\label{test focus} 
\lim_{p\to \infty}(\omega_{tt}^{\hat{ f}_p}(T_1), 
 \omega^{\hat h}(T_1))  = \langle \lambda,\, \varphi(y) \rangle _y .
\eeq
By Theorem \ref{blacho}, the inner products on the 
left side of equation (\ref{test focus}) are obtainable 
from the boundary data. 
Thus, we can find the map $y \to \langle \lambda, \varphi(y) \rangle _y,
\, y \in M^{{\rm int}}$. Since $\varphi \in \Omega^1M^{{\rm int}}$ 
is arbitrary, 
this determines whether $\lambda \in \Omega^1 M^{{\rm int}}$. It also
 determines
whether $\lambda_y =0$ or not for any $y \in M^{{\rm int}}$. 
This  yields  the claim. 
\hfill$\Box$ 
 
Another way to look at Lemma \ref{lem: 9.3} is that 
the admittance map,
${\mathcal Z}^{T_2}$
  determines, for any boundary source
$h \in C_0^{\infty}(]0,T_0[; \Omega^1 \p M)$, the values, at any $y \in M^{{\rm int}}$,
of $\langle \lambda, (\omega ^h_{tt})^1(t)\rangle _y$ for some unknown
$\lambda \in \Omega^1 M^{{\rm int}}$
and $T_1 <t<T_2- {\rm diam}(M)$.  Moreover, using this map,
we can verify that the $1-$forms $\lambda_k (y)$,
corresponding to three families of focusing sources $\tilde{f}_k(y), \, k=1,2,3$,
are linearly independent at any $y \in M^{{\rm int}}$.
These give rise to the following result. 

\begin{lemma}\label{gauge2}
 Let ${\mathcal Z}^{T_2}$ be the admittance map
of the Riemannian manifold with impedance $(M,g,\alpha)$. Then, for 
$T_1\leq t \leq T_2 -{\rm diam}(M)$ and $\hh\in \overline {\F (]0,T_0[)},$
it is possible to find the forms
\ba
L(y)(\omega^\hh_t(y,t))^1,\quad
K(y)(\omega^\hh_t(y,t))^2,
\ea
at any  $ y\in M^{{\rm int}}$.
Here $L(y): T_y^*M\to T_y^*M$ and 
 $K(y): \wedge^2 T_y^*M\to \wedge^2T_y^*M$ are smooth
sections of
${\rm End}(T^*M^{{\rm int}})$ and ${\rm End}(\Lambda^2T^*M^{{\rm int}})$,
correspondingly. 
\end{lemma}

We emphasize that, at this stage, 
$L(y)$  and $K(y)$ are  unknown. However, they are 
independent of $t$ or $\hh$.

{\em Proof:}
As $M$ is already found, we can choose three
 differential $1-$forms, $\xi_k \in 
\Omega^1M^{{\rm int}}$ which, at any $y \in M^{{\rm int}}$, form a basis in $T^*_yM$.
Using the families  $\tilde{f}_k(y)$ of focusing sources introduced earlier,
we can construct, for any $\hh \in \F^\infty(]0,T_0[)$,
the differential $1-$ form,
\beq
\label{25.11.1}
\rho^{\hat{h}}(y,t) := \langle \lambda_k(y), \,
 (\omega_t^{\hat{h}}(y,t))^1 \rangle_y
\, \xi_k(y).
\eeq
This defines a smooth section, $L(y)$ of ${\rm End}(T^*M^{{\rm int}})$,
\ba
L(y)(\omega^\hh_t(y,t))^1 = \rho^{\hat{h}}(y,t)  \in \Omega^1M^{{\rm int}},
\ea
proving  
the assertion for $(\omega^\hh_t(y,t))^1$ with $\hh\in \F^\infty(]0,T_0[)$.
 Its extension to
$\,\hh \in \overline{\F(]0,T_0[)}$ is an immediate corollary of the fact that
$\F^\infty(]0,T_0[)$ is dense in $\overline{\F(]0,T_0[)}$ in the $\F-$norm.

To analyse $(\omega^\hh_t(y,t))^2$ , consider the form 
\beq
\label{dual}
\eta=(\frac 1\alpha *\omega^3,\frac 1\alpha *\omega^2,\frac 1\alpha *\omega^1,
\frac 1\alpha*\omega^0)
\eeq
(cf. $\nu$-forms in formulae (\ref{hodge}) and (\ref {nu field 2})).
This form satisfies the complete Maxwell system
\[
\eta_t+\widetilde{{\mathcal M}} \eta=0,
\]
 where $\widetilde{{\mathcal M}} $ is the differential expression 
(\ref{M}), corresponding to the  manifold
 $(M,g,\alpha^{-1})$.
Then the admittance map $\bt \eta^1|_{\p M\times ]0,T_2[}\to 
\bn_{\alpha^{-1}} \eta^2|_{\p M\times ]0,T_2[}$ is the inverse of
the given admittance 
map ${\mathcal Z}^{T_2}:\bt \omega^1|_{\p M\times ]0,T_2[}\to 
\bn_{\alpha} \omega^2|_{\p M\times ]0,T_2[}$. Thus, ${\mathcal Z}^{T_2}$
determine the admittance map, $\widetilde{{\mathcal Z}}^{T_2}$ for 
 $\widetilde{{\mathcal M}} $ and we can apply the results for 
$(\omega^{\hh}_t)^1$ to
 $(\eta^{\hf}_t)^1$, where $\hf={\mathcal Z}^{T_2}\hh$. 
Namely, we can find $\tilde L(y)(\eta_t(y,t))^1$, where
$L(y): T_y^*M\to T_y^*M$ is a smooth
section of
${\rm End}(T^*M^{{\rm int}})$ which, at this stage, is unknown. 
At last, since 
\ba
*\tilde L(y)(\eta^{\hat f}_t(y,t))^1=K(y)(\omega^{\hat h}_t(y,t))^2,
\ea 
for some smooth section, $K(y)$ of  ${\rm End}(\Lambda^2T^*M^{{\rm int}})$,
the assertion follows. 
\hfill$\Box$ 

We note for the further reference, that, similar to the case of the $1-$forms, 
the construction of $K(y)$ involves a choice of three differential $2-$forms,
that we denote by
$\mu_k 
\in \Omega^2M^{{\rm int}}$, and three families of generalized sources
$\tilde{\kappa}_k(y)$ that satisfy
\beq
\label{25.11.2}
 \omega^{\partial_t  
\tilde{ \kappa_k}}_t(T_1)= (0,\delta(\mu_k\ud_y(x)),0,0),  \quad \mu_k \in \Lambda^2T^*_yM,\quad k=1,2,3.
\eeq
Below we call the generalized sources $\tilde{f}_k(y)$ the 
focusing sources for $2$--forms and
$\tilde{\kappa}_k(y)$ the 
focusing sources for $1$--forms.  

\smallskip

 Before going to a detailed discussion in the next sections of the
reconstruction of $\alpha$, let us explain briefly the main outline of this construction.
It follows from Lemma \ref{gauge2} that, using the admittance map ${\mathcal Z}$, 
we can find the electromagnetic waves $\omega^f_t(t),$ $T_1<t<T_2-\diam(M)$,
  up to  unknown 
linear transformations, $L$ and $K$. We observe that, by Theorem
\ref{global control th}, for any basis $\xi_k(y), \,k=1,2,3$, there are families 
$\tilde{f}_k(y)$ of focusing sources, such that the corresponding transformation $L$ 
is just identity. Indeed, to achieve this goal, we should choose $\tilde{f}_k(y)$ 
in such a manner
that, at any $y \in M^{{\rm int}}$, 
\beq
\label{26.11.1}
(\omega_{tt}^{\tilde{f}_k(y)})^2(T_1) = d(\lambda_k(y) \ud_y),
\eeq
where $\lambda_k(y)$ is dual to $\xi_k(y)$,
\beq
\label{26.11.2}
\langle \lambda_k(y),\, \xi_j(y) \rangle_y = \delta _{kj}.
\eeq
In the next sections we will identify conditions
on $\tilde{f}_k(y)$ and
$\tilde{\kappa}_k(y)$,
verifiable  in terms of ${\mathcal Z}$, which make $L$ and $K$
 to be identities.

\subsection{Reconstruction of the wave impedance} 
 
In the previous section, it was shown how to select 
a family $\tilde{f}_y, \, y \in M^{{\rm int}}$ of focusing sequences 
such that the corresponding electromagnetic  fields 
 concentrate at $t=T_1$
at a single point,  $y$,
\[ 
 \lim_{p\to\infty}\omega_{tt}^{f_p}(T_1)  
 =(0,0,d(\lambda\ud_y),0). 
\] 
Here $\lambda=\lambda_y \in T^*_yM$ is yet unknown. Moreover,
we can select this family in such a manner that, as a function of
$y, \, \lambda_y \in \Omega^1M^{{\rm int}}$. 
In particular, it means that for times 
$0<t<t_y=\tau(y,\partial M)$, the 
electromagnetic
wave defined as 
\begin{equation}\label{el green} 
 G_{\rm e}(y)=G_{\rm e}(x,y,t) = G_{\rm e}[\lambda](x,y,t) 
 = \lim_{p\to\infty}\omega_{tt}^{f_p}(t + T_1), 
\end{equation} 
 satisfies the 
 initial-boundary value problem 
\begin{eqnarray}\label{green equations} 
 (\partial_t + {\mathcal M})G_{\rm e} (y)&=& 0\mbox{ in $M\times]0,t_y[$}, 
\nonumber 
\\ 
\noalign{\vskip4pt} 
 \bt G_{\rm e}(y) &=& 0\mbox{ in $\partial M\times]0,t_y[$},
\\ 
\noalign{\vskip4pt} 
G_{\rm e}(y)|_{t=0} &=& (0,0,d(\lambda\ud_y),0).\nonumber 
\end{eqnarray} 
The solution to this problem is called the {\em electric 
Green's function}. 
We will use this solution  to reconstruct 
the scalar wave impedance, $\alpha$ on $M$. We start
with analysis 
of some properties of  $G_{\rm e}$. To this end, 
we  will represent $G_{\rm e}$ in terms of the 
standard  Green's function, $G=G(x,y,t)$
 for the wave equation
on $1-$ forms.  Thus, $G$ is defined as the solution to the following
initial-boundary value problem
\begin{eqnarray}\label{1-form Green's function} 
(\partial_t^2+d\delta+\delta d)G(x,y,t)&=
& (\partial_t^2+\Delta^1_{\alpha}) G(x,y,t) =\lambda 
\ud_y(x)\ud(t) 
\hbox{ in  $M\times \R$}\nonumber \\ 
\noalign{\vskip4pt} 
 \bt G(x,y,t)=& 0,&
\\ 
\noalign{\vskip4pt} 
G(x,y,t)|_{t<0 }
&=&0,
\nonumber 
\end{eqnarray} 
where $\lambda\in T^*_yM$ is a given 1--form. 
 
This 
Green's function 
has the following asymptotic behaviour. 
 
\begin{lemma}\label{asymptotics of green} 
For $0<t<t_y$,
Green's function $G(x,y,t)$ for the 1--form 
wave equation (\ref{1-form Green's function}), has  
the representation 
\[ 
G(x,y,t)=\ud(t-\tau(x,y))Q(x,y)\lambda+r(x,y,t). 
\] 
Here $Q(x,y):T^*_yM\to T^*_xM$ is a bijective
 map that corresponds to 
a $(1,1)$--tensor depending smoothly on 
$(x,y)\in M^{{\rm int}}\times M^{{\rm int}}\setminus 
\hbox{diag}(M^{{\rm int}})$. The remainder $r(x,y,t)$ is a bounded function, 
when $t<t_y$, where $t_y$ is  small enough. 
\end{lemma} 
 
The proof of this lemma is postponed in the Appendix. 
 
In the following, 
{we  fix $y\in M^{{\rm int}}$ and $\lambda=\lambda_y \in T^*_yM$. 
 By operating with the exterior derivative $d$ on the both 
sides of the differential equation in  
(\ref{1-form Green's function}), we see that 
\[ 
(\partial_t^2+\Delta ^2_{\alpha})dG(x,y,t)=d(\lambda\ud_y) 
\ud(t). 
\] 
Hence, using the decomposition 
$ 
 (\partial_t^2 + {\bf \Delta})= 
 (\partial_t +{\mathcal M})(\partial_t -{\mathcal M}), 
$ 
we find that the form  $\omega(t)=(0,0,dG(x,y,t),0)$ 
satisfies the equation 
\[ 
(\partial_t+{\mathcal M})\bigg((\partial_t-{\mathcal M}) 
(0,0,dG(x,y,t),0)\bigg)=D_{y,\lambda} 
\ud(t), 
\] 
where 
\[ 
 D_{y,\lambda}=(0,0,d(\lambda\ud_y),0). 
\] 
Let $\widetilde{G}_{{\rm e}}(y)= \widetilde{G}_{{\rm e}}(x,y,t)$ 
be defined as 
\begin{eqnarray*} 
 \tilde{G}_{\rm e}(x,y,t) &=&(\partial_t-{\mathcal M}) 
(0,0,dG(x,y,t),0)\\ 
\noalign{\vskip4pt} 
 &=& (0,\delta d G(x,y,t),\partial_t d G(x,y,t),0). 
\end{eqnarray*} 
Then, due to the finite propagation speed, $G(y) \big|_{\p M \times ]0,t_y[} =0$,
so that $\tilde{G}_{\rm e}(x,y,t)$ satisfies the boundary condition 
${\bf t}\tilde{G}_{\rm e}(y) =0$ for $t<t_y$. Invoking the uniqueness of solution for
(\ref{green equations}), we see that $\tilde{G}_{\rm e}(y)=G_{\rm e}(y), \,
t<t_y$.


Now, using Lemma \ref{asymptotics of green}, we obtain that
\ba
dG^1_e(x,y,t)=(Q(x,y)\lambda_y\wedge d\tau(x,y)) \underline \delta^{(1)}
(t-\tau(x,y))+r_1(x,y,t),
\ea
where $\underline{\delta}^{(1)}$ is derivative of the delta-distribution and 
the residual $r_1$ is sum of the 
delta-distribution on $\partial B_y(t)$,
where  $B_y(t)$ is the ball of radius $t$ centered in $y$, 
 and a bounded function. 
 Thus, we see that
\beq
\label{21.11.1}
\nonumber
\p_tdG^1_{{\rm e}}(x,y,t)&=&(Q(x,y)\lambda_y\wedge d\tau(x,y)) \underline \delta^{(2)}
(t-\tau(x,y))+r_2(x,y,t),\\
\nonumber
\delta dG^1_{{\rm e}}(x,y,t)&=&*\left(d \tau(x,y) \wedge
*(Q(x,y)\lambda_y\wedge d\tau(x,y)) \right) \underline \delta^{(2)} 
(t-\tau(x,y)) \\
& &+r_3(x,y,t),
\eeq
where residuals  
$r_2$ and $r_3$  are sums of first and zeroth
derivatives of the delta-distribution on
 $\partial B_y(t)$ and a bounded function. 
Moreover, 
by formulae (\ref{21.11.1}),
\beq
\label{21.11.2}
{\bf t}_{B_y(t)} [ Q(x,y)\lambda_y \wedge d\tau(x,y)]=0,
\eeq
\beq
\label{21.1.3}
{\bf n}_{B_y(t)} \left[*\left(d \tau(x,y) \wedge
*(Q(x,y)\lambda_y\wedge d\tau(x,y)) \right) \right]=0,
\eeq
 where ${\bf t}_{B_y(t)} \omega^k,\, {\bf n}_{B_y(t)}\omega^k$ are the tangential
and normal components of $\omega^k$ on $\p B_y(t)$. 
This corresponds to the physical fact that the wavefronts
of the electric and magnetic fields are perpendicular
to the propagation direction.
Let now  $\tilde{f}_y$ be a focusing source for a point $y \in M^{{\rm int}}$.
Due to the definition (\ref{limit})  and the definition of the generalized source
(\ref{gener}), there is a sequence $(f^y_p)_{p=1}^\infty$, 
$f^y_p \in C^{\infty}_0(]0,T_0[;\Omega^1 \p M),\,
 p=1,2,\cdots$ such that}
\ba
\omega^{\tilde f_y}_t(T_1)=\lim_{p\to \infty}
\omega^{f^y_p}_t(T_1),
\ea
and the right side is understood in the sense of the distribution-forms on 
${\bf \Omega} M^{{\rm int}}$. Then, for $t \geq 0$,
\beq
\label{21.11.4}
\omega^{\tilde f_y}_t(t+T_1)=\lim_{p\to \infty}
\omega^{f^y_p}_t(t+T_1).
\eeq
Applying Lemma \ref{gauge2}, it is possible to find, via given ${\mathcal Z}^{T_2}$,
the magnetic components $K(x)(\omega^{f}_t(x,t+T_1))^2$, 
$f= f^y_p$ of these fields
with $K$ being a smooth section of ${\rm End}(\Lambda^2T^*M^{{\rm int}})$.
At last, using (\ref{21.11.4}), we find
\ba
K(x)(\omega^{\tilde f_y}_t(x,t+T_1))^2 =\lim_{p\to \infty}K(x)(\omega^{ f^y_p}_t(x,t+T_1))^2 \,
\in {\cal D}'(\Omega^2M^{{\rm int}}).
\ea

Since 
\ba
\omega^{\tilde f_y}_{tt}(T_1)=(0,0,d(\lambda \underline\delta_y),0),
\ea
we see that
\ba
\omega^{\tilde f_y}_{tt}(t+T_1)=G_{e}(\cdotp,y,t),
\ea
when $t<t_y$.
In particular, we can find the singularities of  Green's function
up to a linear transformation $K(x)$.

Hence we have shown:
\begin{lemma}\label{values of Maxwell}
Let $\tilde f_y=(f^y_p), \,p=1,2,\cdots,$
be a focusing source for a point $y$. 
Then, 
given the admittance map, ${\mathcal Z}^{T_2}$, it is
possible to find the distribution
$2-$form $K G_{{\rm e}}(y)^2$
for all $x$ satisfying 
$\tau(x,y)<\hat t_y,$ where $\hat t_y$ is small enough.
In particular, the leading singularity of this form
determine the $2-$form 
\beq\label{wave front value}
 K(x)(Q(x,y)\lambda_y\wedge d\tau(x,y)).
\eeq
\end{lemma}

\begin{figure}[htbp]
\begin{center}
\psfrag{1}{$\vec v$}
\psfrag{2}{$\vec w$}
\includegraphics[width=10cm]{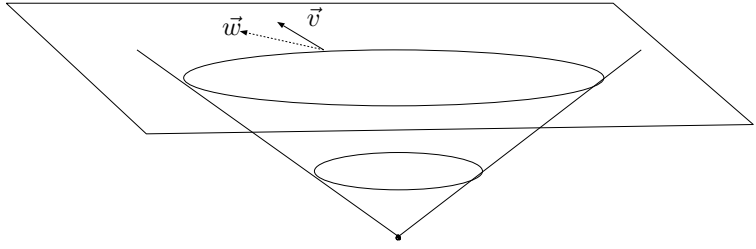} \label{pic 6}
\end{center}
\caption{Vector $\vec v$ is the right singularity of the electromagnetic wave 
in the plane
$M\times \{t\}$. The reconstructed singularity $\vec w$ has
wrong direction, if the transformation matrix $K(x)$ is not isotropic.
}
\end{figure}

 As shown at the end of the previous section, 
$K \in {\rm End}(\Lambda^2T^*M^{{\rm int}})$ was 
obtained by using three focusing sources
$\tilde{\kappa}_{k}(x),$ $k=1,2,3$.
 Our next goal is to formulate conditions, verifiable using 
boundary data, for $K$ to be isotropic, i.e.
\beq\label{K-requirement}
K(x)=c(x)I,
\eeq
with $c(x)$ being a smooth scalar function.

We start with observation that, for a given $\tilde{\kappa}_k(x), \, k=1,2,3,$
and any $\tilde{f}_y$, we can find 
${\bf t}_{B_y(t)}K(\cdot)(Q(\cdot,y)\lambda\wedge d\tau(\cdot,y))\big|_x$,
 for $x \in \partial B_y(t)$ and small $t>0$.  
Here $K(x)$ is the linear transformation corresponding to the chosen
$\tilde{\kappa}_k(x)$. This follows from Lemma 
\ref{values of Maxwell}
and the fact that the underlying Riemannian manifold $(M,g)$ is already found.
When $K$ is isotropic, it follows from (\ref{21.11.2}) that
\beq
\label{orthog}
{\bf t}_{B_y(t)}K(\cdot)(Q(\cdot,y) \lambda \wedge d\tau(\cdot,y))=0.
\eeq
Let us show that the condition (\ref{orthog}), which is verifiable via boundary
data, actually guarantees that $K$ is of form (\ref{K-requirement}). 
}

{Indeed, for a given $\lambda \in T^*_yM$ and any $x,t$ with 
sufficiently small $t=\tau(x,y)$, (\ref{orthog}) means that 
$K(x)(Q(x,y) \lambda \wedge d\tau(x,y))$ continues to be normal
to the $2-$dimensional subspace, $T_x(\partial B_y(t)) \subset T_xM$, i.e.
\ba 
\left[K(x)(Q(x,y) \lambda \wedge d\tau(x,y)) \right](X,Y) =0, 
\quad X,Y \in (T_x\partial B_y(t)).
\ea
As $Q(x,y): T^*_yM \to T^*_xM$ is bijective,
$Q(x,y) \lambda \wedge d\tau(x,y)$ runs through the whole 
$2-$dimensional subspace in $\Lambda^2 T^*_xM$, normal to
$T_x(\partial B_y(t))$, when $\lambda$ runs through $T^*_yM$.
Therefore, (\ref{orthog}) implies that $K(x)$  keeps this
subspace invariant.

Let us now vary $y$ and $t$ keeping $x$ fixed. Then $T_x(\partial B_y(t))$
runs through the whole Grassmannian manifold, $G_{3,2}(T_xM)$.
Therefore, (\ref{orthog}) implies that, if  $\omega^2 \in
\Lambda^2 T^*_xM$ is normal to a subspace 
${\mathcal L} \in G_{3,2}$, then $K(x) \omega^2 $ remains to be normal to 
${\mathcal L}$.
This is implies that the eigenspace of $K(x)$ has dimension three
and hence we must have
\ba
K(x)=c(x)I,
\ea
where $c(x)$ is a scalar function and $I$ is 
the identity map  in $\Lambda^2 T^*_xM$. 

 In the following, we always take  focusing sources which satisfy
 (\ref{K-requirement}). 

Thus, we can find the values of the 2-forms
\ba
(\tilde \omega_t^{\hat{f}}(x,T_1)^2=_{def}c(x) (\omega_t^{\hat{f}}(x,T_1))^2,
\quad \hat f \in {\cal F}.
\ea

Our further considerations are based on the equation, 
\beq
\label{25.11.3}
d(\omega_t^{\hat{f}}(T_1)^2)=0.
\eeq
Thus, we intend to choose those focusing sources $\tilde{\kappa}_k(x)$
which produce $K(x) =c(x)I$ with 
\ba
d( \tilde \omega_t^{\hat{f}}(T_1))^2=0, \quad 
\hat f\in {\cal F}^\infty.
\ea 
In this case,
\ba
0=d(c(x)  (\omega_t^{\hat{f}}T_1))^2)=
dc(x) \wedge  (\omega_t^{\hat{f}}(T_1))^2,
\ea
due to (\ref{25.11.3}). 
By the global controllability, the form  $(\omega_t^{\hat{f}})^2(x,T_1)$ 
runs through the whole $\Lambda^2T^*_xM$. 
 This implies
that   $dc(x)=0$
at any $x \in M^{{\rm int}}$. 
Hence,
$
c(x)=c_0
$
 is a constant. 
Thus, we choose focusing sources, $\tilde{\kappa}_k(x)$, so that 
\ba
(\tilde \omega_t^{\hat h}(T_1))^2=c_0 (\omega_t^{\hat h}(T_1))^2.
\ea
Evaluating the inner products, 
\ba
\int_M \tilde \omega_t^2(x,T_1)\wedge * \tilde \omega_t^2(x,T_1)=
\int_M
c_0 \omega_t^2(x,T_1)\wedge *c_0
\omega_t^2(x,T_1),
\ea
we can compare them with the energy integrals,
\ba
\int_M\frac 1{\alpha(x)} \omega_t^2(x,T_1)\wedge *  \omega_t^2(x,T_1),
\ea
which can be found from the boundary data by means of Theorem
\ref{blacho} . 
By considering waves  $\omega_t^{\hat f_j}(x,T_1)$ with
\ba
\lim_{j\to \infty}
 \hbox{supp}(\omega_t^{\hat f^j}(\cdot,T_1))=\{y\},
\ea 
we  find the ratio
\ba
\lim_{j\to \infty}\frac 
{\int_M (\tilde \omega_t^{\hat f_j}(x,T_1))^2\wedge * 
(\tilde \omega_t^{\hat f_j}(x,T_1))^2}
{\int_M \frac 1{\alpha(x)} (\omega_t^{\hat f_j}(x,T_1))^2
\wedge *  (\omega_t^{\hat f_j}(x,T_1))^2}=
c_0^2\alpha(y).
\ea
The above considerations  imply that, using ${\mathcal Z}^{T_2}$, 
 we can determine $\alpha(x)$ up to a constant
$c_0^2$. Since the impedance map satisfies ${\mathcal Z}_{M,g,c\alpha}=
c^{-1}{\mathcal Z}_{M,g,c\alpha}$, we see that, knowing the impedance map,
we can also determine $c_0$. 

As the fact that ${\mathcal Z}^T, \, T \geq 4\,{\rm diam}M$ 
determines $(M,g)$ is already proven, this completes 
the proof of Theorem \ref{ip}.
\hfill$\Box$ 

\subsection{Back to $\R^3$}
 
In this section we use the obtained uniqueness result
for Maxwell equation on a 3-dimensional manifold to analyze to the 
group of transformations which preserve the boundary data  in the
dynamical inverse problem for Maxwell's equation
(\ref{MF vector}) in a domain $\Omega\subset \R^3$.

Assume that two Maxwell systems 
with electric and magnetic permittivities
$\epsilon^j_k(x)$, $\mu^j_k(x)$ and 
$\tilde \epsilon^j_k(x)$, $\tilde \mu^j_k(x)$, $x\in \Omega$,
have the same admittance map ${\mathcal Z}^T$ on $\p \Omega\times [0,T]$,
where $T$ is sufficiently large.  Denote by $(M,g,\alpha)$
and  $(\tilde M,\tilde g,\tilde \alpha)$ the corresponding
abstract Riemannian manifolds with impedance. By
Theorem \ref{ip}, 
\ba
(\tilde M,\tilde g,\tilde \alpha)=(M,g,\alpha)
\ea
i.e., there is an isometry $H:(M,g)\to (\tilde M,\tilde g)$
and $\alpha=H_*\tilde \alpha$. 
We can represent  the abstract manifold 
$(M,g,\alpha)$ as the domain $\Omega$ with the metric tensor, $g_{ij}(x)$,
given  in
Euclidean coordinates by (\ref{connection e and g}), 
(\ref{travel time metric}),
and  scalar impedance, $\alpha(x)$  in these coordinates. 
Similarly, we  represent manifold 
$(\tilde M,\tilde g,\tilde \alpha)$ using Euclidean coordinates in $\Omega$
and obtain the metric tensor
$\tilde g_{ij}$ and impedance $\tilde \alpha$. Then 
$H:\tilde M\to M$ corresponds to a diffeomorphism,
\beq\label{19.0}
\tilde X:\Omega\to \Omega,\quad \tilde X|_{\p \Omega}=id|_{\p \Omega},
\eeq
and
\beq\label{19.1}
& &\tilde g=\tilde X_* g,\quad\hbox{i.e., } \,\,
\tilde g^{ij}(\tilde x)=\frac {\p\tilde x^i} {\p x^p}
\frac {\p\tilde x^j} {\p x^q}g^{pq}(x),\quad \tilde x=\tilde X(x),\\ 
\nonumber
& &\tilde \alpha=\tilde X_* \alpha,\quad\hbox{i.e., } \,\,
\tilde \alpha(\tilde x)=\alpha(x).
\eeq
Using (\ref{19.1}) and (\ref{travel time metric}),
we see also that 
\beq\label{19.2}
& &\tilde g_\epsilon=\tilde X_* g_\epsilon,\quad\hbox{i.e., } \,\,
\tilde g^{ij}_\epsilon(\tilde x)=\frac {\p\tilde x^i} {\p x^p}
\frac {\p\tilde x^j} {\p x^q}g^{pq}_\epsilon(x),\quad \tilde x=\tilde X(x). 
\eeq
Employing  formula (\ref{connection e and g}),
we obtain
\beq\label{19.2b}
\epsilon^p_q=\sqrt{g_\epsilon} g_\epsilon^{pr}\delta_{pq},\quad
\tilde \epsilon^p_q=\sqrt{\tilde g_\epsilon} \tilde g_\epsilon^{pr}\delta_{pq}.
\eeq
Combining formulae (\ref{19.1})--(\ref{19.2b})
and introducing
\ba
\epsilon^{pq}=\epsilon^p_r\delta^{jq},\quad
\tilde \epsilon^{pq}=\tilde \epsilon^p_r\delta^{jq},
\ea
we obtain
\beq\label{19.3}
\tilde \epsilon^{pq}=\frac 1{\hbox{Det}\,(D\tilde X)}
\frac {\p\tilde x^i} {\p x^p}
\frac {\p\tilde x^j} {\p x^q}\epsilon^{pq}(x),\quad \tilde x=\tilde X(x).
\eeq
Similarly,
\beq\label{19.4}
\tilde \mu^{pq}=\frac 1{\hbox{Det}\,(D\tilde X)}
\frac {\p\tilde x^i} {\p x^p}
\frac {\p\tilde x^j} {\p x^q}\mu^{pq}(x),\quad \tilde x=\tilde X(x)
\eeq
Clearly, if $\tilde X:\Omega\to \Omega$ and $\tilde X|_{\p \Omega}=id$,
the admittance map ${\mathcal Z}^T,$ $T>0$ is preserved in
transformations (\ref{19.0}),  (\ref{19.3}),  (\ref{19.4}).

Thus we have proven the following result.

\begin{theorem} 
\label{group}
The group of transformations
for Maxwell's equations (\ref{MF vector})--(\ref{constitutive}) 
with a scalar wave impedance, 
which preserves
the admittance map  ${\mathcal Z}^T,\,$ $T>4\,{\rm diam}(M,g)$,
is generated by the group of diffeomorphisms of $\Omega$ satisfying
(\ref{19.0}). The corresponding transformations of $\epsilon$ and $\mu$
are then
defined by formulae
 (\ref{19.3}),  (\ref{19.4}). 
\end{theorem}

{\bf Remark 9.} It follows from (\ref{19.3}),  (\ref{19.4}), that 
$\epsilon^{jk}$ and  $\mu^{jk}$ 
 do not transform like tensors. This
is due to the special role played by the underlying Euclidean
metric $g_0^{ij}=\delta^{ij}$, which does not change
by diffeomorphisms $\tilde X$. It should be noted that
this form of transformations is observed also in the study
of the Calder\'{o}n inverse conductivity problem.
Indeed, it is shown in \cite{Sy} that, for the  conductivity
equation in $\Omega\subset \R^2$, the boundary measurements
determine the anisotropic conductivity up to same group
of transformations 
as described in Theorem \ref{group}. 
The Calderon problem
is closely related to the inverse problem for Maxwell's equation,
for instance, the low-frequency limit of the admittance
map ${\mathcal Z}^\infty$ is related to the Dirichlet-to-Neumann
map for the conductivity equation \cite{La3}.

\subsection{Outlook}

There are several direction to which the present work can 
be extended. 

1. Natural question is the minimal observation time required
to parameter reconstruction. It can be shown that the admittance
map ${\mathcal Z}^t$ for any $t>0$.  
Thus, it follows from the above,
 that  ${\mathcal Z}^T$, $T>2\,{\rm rad}M$ determines uniquely
the manifold $M$, 
 metric $g$ and wave impedance $\alpha$. 
The reconstruction of  ${\mathcal Z}^t$ for any $t>0$
may be  obtained
by a direct continuation of the admittance map, i.e.,
without solving the inverse problem. 
This continuation
is 
a direct generalization to the considered Maxwell's case of
 the technique
developed in \cite{KL1}, \cite{KKL} 
for the scalar wave equation.
An analogous method has recently
be applied to Maxwell's equations in \cite{BI}.

2. Another natural inverse boundary value problem is the
inverse boundary spectral problem for the electric
Maxwell operator $\M_{\rm e}$ defined in Definition 
\ref{d. 2}. The problem is to determine the metric
$g$ and wave impedance $\alpha$, or, in the other words, $\e$ and $\mu$
from the non-zero eigenvalues $\lambda_j$ of $\M_{\rm e}$
and the normal boundary values of the corresponding 
eigenforms. This problem was studied  in, e.g. \cite{L1}, \cite{L2},
for scalar
Maxwell's equations. 
For 
the considered anisotropic case,
 this requires significant 
modifications of the method developed in this paper
and will be published elsewhere.

3. It often occurs  in applications 
that the measurements
are made only on a part of boundary. In formalism of this
paper this means 
that we actually know only the restriction of the admittance
map to the part $\Gamma\times [0,T]$ of the lateral boundary, i.e.
we are given
\ba
 {\mathcal Z}^Tf|_{\Gamma\times [0,T]}, \quad
f \in C^{\infty}_0([0,T]; \Omega^1\Gamma). 
\ea
For the scalar wave
equation the corresponding problem is studied in \cite{KK1}, \cite{KK2}
(see also  \cite{KL3},  \cite{KKL}, \cite{KKL2}).
The combination of these methods 
and those of the present paper will be useful for analyzing
the corresponding problem for Maxwell's equations.

\section*{Appendix: The WKB approximation} 

Here we consider asymptotic results 
for Green's function
or, more precisely, for  Green's 1-form,
$G(x,y,t)=G_{\lambda}(x,y,t)$,
which is defined as the solution 
for the wave equation
\beq
\label{11.01}
& &\p_t^2 G_{\lambda}  + (d \delta  + \delta d) G_{\lambda} =
a  \ud _y (x)\ud  (t)
\quad \hbox{in} \,\, M\times \R_+,\\
& & \nonumber
G_{\lambda} (\cdot,y,t) = 0 \quad \hbox{for} \,\, t<0, \quad  {\bf t}G_{\lambda} (\cdot,y,t) = 0,
\eeq
where $\delta = \delta_{\alpha}$.
Here, $\lambda $ is a 1-form
$
\lambda = \sum_{i=1}^3 \lambda_i \, dx^i
$
in normal coordinates $(B_y(\rho),X),\,X=(x^1,x^2,x^3)$,  near 
a point $y \in M^{{\rm int}}, \, X(y)=0$.
We assume that 
$B_y(\rho) \cap \partial M = \emptyset$. 
 In these coordinates, 
 $ \ud _y (x) =  \ud (x)$,
when $x \in U$. Clearly, we can find, instead of the solution to
(\ref{11.01}), the fundamental solution
\beq
\label{11.06}
& &\p_t^2 G + (d \delta  + \delta  d) G =
I  \ud (x) \ud (t)
\quad \hbox{in} \,\, M\times \R_+,\\
& & \nonumber
G(\cdot,y,t) = 0 \quad \hbox{for} \,\, t<0, \quad {\bf t}G(\cdot,y,t) = 0,
\eeq
where $I$ is the $3 \times 3$ identity matrix.
Equation (\ref{11.06}), written in normal coordinates,
becomes a hyperbolic system
\beq
\label{11.02}
\left\{ (\p_t^2 - g^{ij} \p_i \p_j)I + B^i \p_i +C \right \}
G = I \ud (x) \ud (t).
\eeq
Here $g^{ij}(x)$ is the metric tensor in these coordinates with
\beq
\label{11.05}
g^{ij}(0)= \delta^{ij}, \, \p_kg^{ij}(0)=0,
\eeq
and $B^i(x), \, C(x)$ are smooth
$3 \times 3$ matrices.
We note that,
 in normal coordinates, $\tau(x,y) =|x|$. However, we prefer to
keep the notation $\tau$ to stress the invariant nature of considerations
 below.

Following \cite{Co}, \cite{Ba},
which deal with the scalar case, 
we search for the solution to (\ref{11.02}) in the WKB form:
\beq
\label{11.03}
G(x,t) = G_0(x) \, \ud (t^2 - \tau^2) +
\sum_{l \geq 1} G_l(x) \, S_{l-1} (t^2 - \tau^2),
\eeq
where
$
S_l(s) = s_+^l /\Gamma (l+1).
$
Substitution of  expression (\ref{11.03}) into equation (\ref{11.02})
gives rise to the recurrent system of (transport) equations.
The principal one
is the equation for $G_0$,
\beq
\label{11.04}
4 \tau \frac{d G_0}{d \tau}(\tau \hat{x}) +
\left \{(g^{ij} (\tau \hat{x}) \, \p_i\p_j \tau ^2(\tau \hat{x}) - 6)\,I
  + B^i(\tau \hat{x}) \, \p_i \tau^2(\tau \hat{x}) \right \}
  G_0(\tau \hat{x}) = 0,
\eeq
where $\hat{x} = x/\tau.$ In addition, to satisfy initial
conditions $I \ud (x) \ud (t)$, corresponding  to the right side in the 
wave equation 
(\ref{11.06}), 
  we require that
\beq
\label{12.02}
G_0(0) = \frac{1}{2 \pi} I.
\eeq
By (\ref{11.05}), $g^{ij} \p_i\p_j \tau ^2 - 6$ is a smooth function
near $x=0$ and $g^{ij} \p_i\p_j \tau ^2\big|_{x=0} - 6=0 $. Also,
$\, \p_i \tau^2 \big|_{x=0} =0$. Therefore,
\bfo
\frac{1}{4 \tau}\left \{(g^{ij} (\tau \hat{x}) \, \p_i\p_j \tau 
^2(\tau \hat{x}) - 6)\,I
  + B^i(\tau \hat{x}) \, \p_i \tau^2 (\tau \hat{x}) \right \}
\efo
is a smooth function of $(\tau, \hat{x})$, so that
  $G_0(x)$ is a smooth $3 \times 3$ matrix of $x$ for $\tau >0$.

{\nottopapertext
\smallskip
\noindent Actually, the matrix $G_0(x)$ is smooth everywhere in the 
neighborhood
of $y$, i.e. $x=0$,  including $x=0$ itself. Indeed, if we write 
the Taylor expansion  of
$(g^{ij}(x) \p_i\p_j \tau ^2(x) - 6) \,I + B^i(x) \p_i \tau^2(x)$ near
$x=0$ and divide the result by  $\tau = |x|$, we obtain that
\beq
\label{12.01}
(g^{ij}(x) \p_i\p_j \tau ^2(x) - 6) \,I + B^i(x) \p_i \tau^2(x) =
\sum_{|\beta| \geq 1} D_{\beta} \tau^{|\beta| - 1}
\hat{x}^{\beta}.
  \eeq
Substituting the Taylor expansion (with respect to $\tau$)
of $G_0(\tau, \hat{x})$,
\bfo
G_0(\tau, \hat{x}) = \sum_{p \geq 0} G_{0,p}(\hat{x})\tau ^p,
\efo
into (\ref{11.04}) and using (\ref{12.01}), (\ref{12.02}), we obtain that
$G_{0,p}(\hat{x})\tau ^p$ are 
homogeneous polynomials of $x$ of degree $p$:
\bfo
G_{0,p}(\hat{x})\tau ^p = \sum_{|\beta| = p} G_{0,\beta}x^{\beta}.
\efo
Then,
$
G_0^p(x) = \sum_{|\beta| < p} G_{0,\beta}\,x^{\beta}
$
satisfies
\beq
\label{12.04}
4 \tau \frac{d G_0^p}{d \tau} +
\left \{(g^{ij}  \, \p_i\p_j \tau ^2 - 6)\,I
  + B^i \, \p_i \tau^2 \right \}
  G_0^p =
\theta^p,
\eeq
where
\beq
\label{12.03}
\theta^p(\tau, \hat{x}) = \theta^p(x) \in C^{\infty}(U),
\quad \theta^p(x) = O(\tau^p).
\eeq
We construct $G_0$ as $G_0^p \,(I+  \tilde{G}_0^p)$. Substituting
this expression into
(\ref{11.04}) and using (\ref{12.04}), (\ref{12.03}), we obtain that
\bfo
4 \tau \frac{d \tilde{G}_0^p}{d \tau}(\tau,\hat{x}) = A^p(\tau,\hat{x}) +
A^p(\tau,\hat{x})\, \tilde{G}_0^p(\tau,\hat{x}),
\quad \tilde{G}_0^p(0)=0,
\efo
where
\bfo
A^p(\tau,\hat{x}) = - \left(G_0^p(\tau,\hat{x})\right)^{-1} \,
\theta^p(\tau,\hat{x}).
\efo
Therefore,  $A^p(\tau,\hat{x}) \in C^{\infty}$ as
a function of $(\tau, \hat{x})$ and
$A^p(\tau) = O(\tau^p)$. 
This implies that $\tilde{G}_0^p(\tau,\hat{x}) = O(\tau^p)$
and is
$C^{\infty}$ smooth as a function of $ (\tau,\hat{x})$, so that
$\tilde{G}_0^p $,
considered as a function of $x =\tau \hat{x}$ is in
$ C^p(B_y(\rho))$.  As $p >0$
is arbitrary and the solution $G_0$ of (\ref{11.04}), (\ref{12.02}) is unique,
$G_0 \in C^{\infty}(U)$.

For
$G_l, \, l \geq 1$, we obtain  transport equations 
\ba
& &4 \tau \frac{d G_l}{d \tau} +
\left \{(4l -6 + g^{ij}(x) \p_i\p_j \tau ^2(x)) \,I
+ B^i(x) \p_i \tau^2 \right \} G_l\\
& & =
  \left[   g^{ij} \p_i \p_jI - B^i \p_i -C \right ]
G_{l-1},
\ea
and $G_l(0)=0$.
If we write $G_l = G_0 F_l$, we obtain for  $F_l$ the equations
\beq
\label{11.07}
4 \tau \frac{d F_l}{d \tau} +4l\,F_l =
G_0^{-1} \, \left[   g^{ij} \p_i \p_jI - B^i \p_i -C \right ]
G_{l-1}, \quad F_l(0) =0,
\eeq
with their solutions
\beq
\label{11.08}
F_l(x) = \frac14 \tau ^{-l} \int _0^{\tau}
G_0^{-1}(s \hat{x}) \, \left\{
\left[ g^{ij} \p_i \p_jI - B^i \p_i -C \right ]
G_{l-1} \right \}(s \hat{x}) \, s ^{l-1} d s,
\eeq
being a smooth function of $x$.

As (\ref{11.02}) is a hyperbolic system, it is easy to 
show that the right side of (\ref{11.03}) represents the asymptotics
 with respect to smoothness of the Green's $1-$form $G(x,y,t)$,
when $t < \tau(y,\partial M)$. 
Clearly, (\ref{11.03})
can also be written in the form
\bfo
G(x,t) = G_0(x) \, \ud (t^2 - \tau^2) +r(x,t)
\efo
where $r(x,t)$ is a bounded
$3 \times 3$ matrix.

\end{document}